\documentclass[11pt,a4paper]{article}

\usepackage{comment} 
\usepackage{graphicx}
\usepackage{amsmath}
\usepackage{amsfonts}
\usepackage{amssymb}
\usepackage{enumerate}
\usepackage{lscape}
\usepackage{longtable}
\usepackage{rotating}
\usepackage{multirow}
\usepackage{color}
\usepackage{url}
\usepackage{subfigure}
\usepackage{rotating}
\newtheorem{theorem}{Theorem}[section]

\usepackage[ruled,vlined,noline,linesnumbered]{algorithm2e}

\usepackage{titlesec}
\usepackage[titletoc]{appendix}

\newtheorem{lemma}{Lemma}[section]

\newtheorem{proposition}{Proposition}[section]

\newtheorem{assumption}{Assumption}[section]

\newenvironment{proof}[1][Proof]{\textbf{#1.} }{\ \rule{0.5em}{0.5em} \vspace{1ex}}

\usepackage{algorithmic}

\def\real{\mathbb{R}}
\def\R{\mathbb{R}}
\def\N{\mathbb{N}}
\def\T{\mathrm{T}}
\def\flow{f_{\mbox{\rm\scriptsize low}}}
\def\calN{\mathcal{N}}
\def\calL{\mathcal{L}}
\def\calS{\mathcal{S}}
\def\calU{\mathcal{U}}
\def\bigx{\mathbf{x}}
\def\bigd{\mathbf{d}}
\def\bigv{\mathbf{v}}

\def\ks{{(k)}}
\def\k1s{{(k+1)}}

\def\Am{A_m}

\DeclareMathOperator{\one}{\mathbf{1}}
\DeclareMathOperator{\eye}{\mathbf{I}}

\usepackage[colorinlistoftodos,bordercolor=orange,backgroundcolor=orange!20,linecolor=orange,textsize=scriptsize]{todonotes}

\DeclareMathOperator{\cm}{cm}

\title{Direct-search methods for decentralized blackbox optimization}
\date{March 5, 2026}
\author{E. Bergou\thanks{Mohammed VI Polytechnic University, Ben Guerir, Morocco.
 ({\tt elhoucine.bergou@um6p.ma}).}
\and Y. Diouane \thanks{GERAD and Department of Mathematics and Industrial Engineering, Polytechnique Monteal, Montreal, Canada. 
(\texttt{youssef.diouane@polymtl.ca}). Funding for this author's research 
was partially provided by the NSERC Discovery grant (RGPIN-2024-0509).}
\and V. Kungurtsev \thanks{Department of Computer Science, Faculty of Electrical Engineering, Czech Technical University in Prague. 
(\texttt{vyacheslav.kungurtsev@fel.cvut.cz}).}
  \and C. W. Royer \thanks{LAMSADE, CNRS, Universit\'e Paris Dauphine-PSL, 
  75016 Paris, France.  
(\texttt{clement.royer@lamsade.dauphine.fr}). 
Funding for this author's research 
was partially provided by CNRS under the IEA grant BONUS and by Agence 
Nationale de la Recherche through program ANR-23-IACL-0008 (PR[AI]RIE-PSAI).}}

\begin{document}
\maketitle

\begin{abstract} 
Derivative-free optimization algorithms are particularly useful for tackling blackbox optimization problems where the objective function arises from complex and expensive procedures that preclude the use of classical gradient-based methods. In contemporary decentralized environments, such functions are defined locally on different computational nodes due to technical or privacy constraints, introducing additional challenges within the optimization process. 

In this paper, we adapt direct-search methods, a classical technique in derivative-free optimization, to the decentralized setting. In contrast with zeroth-order algorithms, our algorithms rely on positive spanning sets to define suitable search directions while still 
possessing global convergence guaranties, thanks to carefully chosen stepsizes. Numerical experiments highlight the advantages of direct-search techniques over gradient-approximation-based strategies.
\end{abstract}

\section{Introduction}
\label{sec:intro}

Decentralized optimization (also referred to as distributed, network, or 
consensus optimization in the literature) has been an increasingly popular topic 
of investigation in recent years~\cite[Chapter 11]{EKRyu_WYin_2022}, 
as networked control and learning systems have proliferated. 
In a decentralized setting, an objective function is defined as a sum of 
functions wherein each individual function is known only locally to some agent, 
and the problem can be solved only through peer-to-peer communication. Note 
that this is distinct from another form of distributed optimization, namely ``federated'' 
optimization, wherein a central node uses subsidiary worker nodes to ease the 
computational burden, but otherwise coordinates the 
procedure~\cite{gao2023decentralized}. 
Decentralized optimization algorithms work at the agent level by updating 
local copies of the problem variables individually, then combining these 
updates with information obtained through communications. The most classical 
algorithm of that form is decentralized gradient descent (DGD), for which a 
number of convergence results have been derived in the convex 
setting~\cite{nedic2009distributed,yuan2016convergence}. The nonconvex case  
proved significantly more challenging, leading to the development of another 
class of algorithms termed gradient tracking~\cite{di2016next,
liu2024decentralized,shah2024adaptive}. Still, convergence guaranties have 
also been derived for DGD techniques in the nonconvex setting, either by 
relying on carefully chosen step sizes~\cite{zeng2018nonconvex} or focusing 
on subclasses of communication networks~\cite{sun2016distributed}. A key 
distinction between the aforementioned results and their centralized 
counterparts lies in the fact that each agent possesses its own copy of the 
problem variables. As a result, a globally convergent method should guaranty 
that all copies are eventually in agreement, i.e., that \emph{consensus} is 
reached among all agents. This property is typically obtained through the 
analysis of an appropriate Lyapunov function, possibly used within the 
algorithm itself~\cite[Chapter 11]{EKRyu_WYin_2022}. An alternate approach 
consists in ensuring approximate consensus by a penalty reformulation of 
the problem, where the level of consensus is controlled by the penalty 
parameter~\cite{zeng2018nonconvex}. In both cases, improvement of the 
objective function is obtained for each agent by taking steps in the 
negative gradient directions.

When derivatives of the local function are unavailable, the decentralized 
optimization literature has focused on using zeroth-order algorithms that 
estimate gradients or directional derivatives through finite-difference 
type estimates~\cite{dang2024adaptive,ghadimi2013stochastic,
hajinezhad2017zeroth,li2021communication,sahu2018distributed,
sahu2020decentralized,tang2020distributed}. Although a common choice to 
alleviate the absence of explicit gradient values, zeroth-order approaches 
are only a subset of derivative-free optimization, a field that has grown in 
importance due to multiple applications in engineering simulations and parameter tuning~\cite{AuHa2017,conn2009introduction,
JLarson_MMenickelly_SMWild_2019}. Among derivative-free optimization 
techniques, direct-search algorithms proceed by exploring the variable 
space through suitably chosen directions. Even though those directions are 
chosen without gradient knowledge, the resulting algorithms can be endowed 
with a rich convergence analysis~\cite{KJDzahini_FRinaldi_CWRoyer_DZeffiro_2024,
kolda2003directsearch}. In addition, direct-search schemes have been 
successfully implemented in parallel environments, with convergence 
guarantees being established in both synchronous and asynchronous 
settings~\cite{JDGriffin_TGKolda_RMLewis_2008,PDHough_TGKolda_VTorczon_2001,
TGKolda_2005}. However, to the best of our knowledge, 
these results only consider a centralized setting, and an investigation (both 
theoretical and practical) of direct-search methods on decentralized 
optimization problems has yet to be conducted.

In this paper, we propose variants on the direct-search paradigm dedicated to 
solving decentralized optimization problems. We preserve key features of 
direct-search schemes, such as the use of positive spanning sets to explore 
the variable space and sufficient decrease conditions, but we adapt stepsize 
conditions in order to guarantee convergence despite the decentralized 
environment. Our analysis either guarantees convergence to at least a stationary 
point of a certain penalty function corresponding to the problem, or certifies 
consensus in the limit. Our numerical experiments, in which we adapt a standard 
DFO benchmark to the decentralized setting, show promising performance of 
direct-search schemes compared to zeroth-order strategies.

The rest of this paper is organized as follows. Section~\ref{sec:background} 
formalizes the decentralized optimization setting, and provides background 
material on decentralized gradient descent and its zeroth-order variants. 
Section~\ref{sec:algos} details two direct-search proposals based on adapting 
the decentralized gradient iteration. Section~\ref{sec:cv} contains global 
convergence results for the two proposed methods. Section~\ref{sec:num} 
reports numerical comparisons between our algorithms and zeroth-order 
strategies. We discuss our findings in Section~\ref{sec:conc}.


\paragraph{Notations}
 In the rest of the paper, $\|x\|$ will denote the Euclidean norm of the 
vector $x$; we will use the same notation for the induced operator norm on 
matrices. For any $r \in \N$, the vector of ones in $\R^r$ will be indicated 
by $\one_r$, while $\eye_r$ will denote the identity matrix in 
$\R^{r \times r}$. Given two matrices $A \in \R^{n_1 \times n_2}$ and 
$B \in \R^{n_3 \times n_4}$, $C:=A \otimes B$ denotes the Kronecker product of 
$A$ and $B$, i.e. the matrix $C \in \R^{n_1 n_3 \times n_2 n_4}$ given by 
$C= \left[ 
\begin{array}{ccc} 
a_{11} B &\cdots & a_{1 n_2} B \\ 
\vdots &\ddots &\vdots \\
a_{n_1 1} B &\cdots &a_{n_1 n_2} B
\end{array}
\right]$.

\section{Background on decentralized optimization}
\label{sec:background}

In this paper, we are interested in the following optimization problem
\begin{equation}\label{eq:pb}
	\min_{x \in \R^n} ~~\sum_{i=1}^m f_i(x),
\end{equation}
where every function $f_i$ is smooth, but its derivative is unavailable for 
algorithmic use. Throughout this paper, we make the following standard
assumptions about the objective function.

\begin{assumption}
\label{as:lip}
For every $i \in \{1,\dots,m\}$, the function $f_i$ is continuously 
differentiable and the gradient $\nabla f_i$ is $L_i$-Lipschitz continuous.
\end{assumption}

\begin{assumption}
\label{as:flow}
The function $x \mapsto \sum_{i=1}^m f_i(x)$ is bounded from below by
$f_{\mathrm{low}} \in \R$.
\end{assumption}

In a decentralized setting, the problem~\eqref{eq:pb} is to be solved over a
network of $m$ agents modeled by an undirected graph 
$\mathcal{G}=(\mathcal{V}, \mathcal{E})$, where 
$\mathcal{V}=\{1, \ldots, m \}$ represents the set of agents and
$\mathcal{E}\subseteq \mathcal{V}\times\mathcal{V}$ represents the set of
edges. Each agent $i$ can evaluate the function $f_i$, but is unaware of
the other functions that define the objective.
In order to solve problem~\eqref{eq:pb} over the network, every agent 
$i \in \mathcal{V}$ thus maintains its own copy of the vector of parameters, 
denoted by $x_i \in \R^n$. 
Letting
$
	\bigx:= \left[ \begin{array}{c} x_1 \\ \vdots \\ x_m \end{array} \right] \in \R^{mn}
$
denote the concatenation of all local copies, the original problem~\eqref{eq:pb} 
becomes equivalent to
\begin{equation}
\label{eq:conspb}
	\begin{array}{ll}
\displaystyle	\min_{\bigx \in \R^{mn}} & F(\bigx) := \sum_{i=1}^m f_i(x_i) \\
	\mbox{subject to} &x_i = x_j \quad \forall (i,j) \in \mathcal{E}.
	\end{array}
\end{equation}

To enforce the constraints in problem~\eqref{eq:conspb}, any agent $i$ can
obtain the values of the copies of its immediate neighbors through
communication. 
Throughout this paper, we assume that the communication 
graph is connected so that information can circulate through the entire 
graph, and meaningful notions of convergence and consensus across agents can 
be defined.
Letting $\calN_i:= \left\{j \in \{1,\dots,m\}
\ \middle|\ j \neq i,\ (i,j) \in \mathcal{E} 
\right\}$ denote the sets of immediate neighbors of agent $i$, a round of 
communications allows agent $i$ to receive
$x_{\calN_i} = \{x_j\}_{j \in \calN_i}$. 
This information is typically 
combined with that of agent $i$ by means of a \emph{mixing matrix} 
$W \in \R^{m \times m}$ whose non-zero entries correspond to elements in $\{\calN_{i} \}$, for which we enforce the following requirements.


\begin{assumption}[Mixing matrix]
\label{as:W}
Let $\lambda_1(W) \ge \lambda_2(W) \ge \dots \ge \lambda_m(W)$ denote the eigenvalues of the mixing matrix $W=[w_{ij}] \in \R^{m \times m}$. 
The matrix $W$ satisfies the following conditions:
\begin{enumerate}[i)]
    \item $W$ is symmetric with nonnegative entries, and $w_{ij} > 0$ if and only if $i=j$ or $j \in \mathcal{N}_i$;
    \item $\lambda_1(W)=1$ and $\lambda_2(W) < 1$;
    \item $W\mathbf{1}=\mathbf{1}$;
    \item $-1 < \lambda_m(W) \leq 0$.
\end{enumerate}
\end{assumption}
When the underlying communication graph is connected, it is sufficient to choose $W$ as a symmetric matrix with nonnegative entries such that 
$w_{ij} > 0$ if and only if $i=j$ or $j \in \mathcal{N}_i$, which is also doubly stochastic, i.e.
\[
	\sum_{i=1}^m w_{ij} = 1 \quad \forall j \in \{1,\dots,m\},
	\qquad
	\sum_{j=1}^m w_{ij} = 1 \quad \forall i \in \{1,\dots,m\}.
\]

Under Assumption~\ref{as:W}, the eigenvalues of $W$ necessarily lie between $-1$ 
and $1$, with the largest eigenvalue equal to $1$. As a result, 
problem~\eqref{eq:conspb} is equivalent to
\begin{equation}
\label{eq:conspbW}
	\begin{array}{ll}
\displaystyle	\min_{\bigx \in \R^{mn}} &F(\bigx)= \sum_{i=1}^m f_i(x_i) \\
	\mbox{subject to} &x_i = \sum_{j=1}^m w_{ij} x_j \quad \forall i=1,\dots,m.
	\end{array}
\end{equation}

It follows that any solution $\bigx^* = [x_i^*]_i \in \real^{nm}$ of 
problem~\eqref{eq:conspb} must satisfy $(\eye_{nm}-\widehat W)\bigx^*=0$, 
where $\widehat W:= W \otimes \eye_n \in \R^{nm \times nm}.$

\subsection{Penalty function and reformulation}
\label{subsec:penalty}

A common approach to decentralized optimization consists of replacing the
constrained formulation~\eqref{eq:conspb} with an unconstrained minimization
problem of a suitable penalty function, where the penalty function is 
typically chosen to be quadratic and depends on the mixing matrix 
$W$~\cite{zeng2018nonconvex}. The resulting optimization problem has the 
form
\begin{equation}\label{eq:penpb}
	\min_{\bigx\in\mathbb{R}^{nm}} ~~
	\calL(\bigx;\gamma) := F(\bigx)+P(\bigx;\gamma), 
\end{equation}
where $F(\bigx):=\sum_{i=1}^m f_i(x_{i})$, $\gamma>0$, and
\[
    P(\bigx;\gamma)
    := \frac{1}{2\gamma}\|\bigx\|^2_{\eye_{nm}-\widehat W} 
    =  \frac{1}{2\gamma}\bigx^\T (\eye_{nm}-\widehat W) \bigx
    =  \frac{1}{2\gamma}\sum_{i=1}^m \|\widehat x_i-x_i\|^2,
\]
with $\widehat x_i := \displaystyle \sum_{j \in \calN_i \cup \{i\}} w_{ij}x_j$. 
The gradient of this penalty function is given by
\[
    \nabla P(\bigx;\gamma)= \frac{1}{\gamma}(\eye_{nm}-\widehat W) \bigx. 
\]

As $\gamma\to 0$, solutions of the penalized formulation~\eqref{eq:penpb} 
converge to those of the constrained problem~\eqref{eq:conspb}. However, 
note that a stationary point of $\bigx^*=[x_i^*]$ in problem~\eqref{eq:penpb} 
is not necessarily a stationary point of problem~\eqref{eq:conspb}, since 
the vectors $x_1^*,\dots,x_m^*$ need not be identical. Still, standard analyzes 
of decentralized gradient methods establish convergence towards a stationary 
point of the penalty function~\cite{nedic2009distributed,zeng2018nonconvex}. 

\subsection{Decentralized gradient descent}
\label{ssec:dgd}

Decentralized gradient descent is based on the following recursion
\begin{equation}
\label{eq:dgd}
    x_{i}^\k1s 
    = \widehat x_i^{(k)}- \alpha^\ks \nabla f_i(x_i^\ks)
    \qquad 
    \forall k \in \N,\ \forall i \in \{1,\dots,m\},
\end{equation}
where $\displaystyle\widehat x_i^{(k)}:= \sum_{j \in \calN_i \cup \{i\}} w_{ij} x_j^\ks $ 
and $\alpha^\ks>0$ is a positive stepsize. In a decentralized setting, the 
stepsize sequence $\{\alpha^{\ks}\}$ is typically fixed a priori, either as 
a constant or a decreasing sequence, and all agents use the same sequence 
throughout the iterations. 

A natural extension of decentralized gradient techniques in the absence of 
derivatives consists of approximating derivatives through deterministic or 
randomized finite differences. The resulting algorithms, termed 
\emph{zeroth-order} decentralized gradient techniques, have the form
\begin{equation}
\label{eq:zodgd}
    x_{i}^{\k1s} 
    = 
    \widehat x_i^{\ks} 
    - \alpha^\ks g_i(x_i^\ks) 
    \qquad 
    \forall k \in \N,\ \forall i \in \{1,\dots,m\},
\end{equation}
where $g_i(x_i^\ks)$ is a gradient approximation. Akin to their first-order 
counterparts, zeroth-order decentralized gradient methods do not use function 
values to check for a decrease in the objective. Although appropriate in a 
decentralized environment, this paradigm differs significantly from the 
dominant approach in derivative-free optimization, which consists of accepting 
steps that reduce the objective value and rejecting those that do 
not~\cite{conn2009introduction}. In the next section, we describe two ways to 
adapt this approach to the decentralized setting.

\section{Decentralized direct-search frameworks}
\label{sec:algos}

In this section, we propose two ways to adapt the classical direct-search 
algorithmic framework~\cite{kolda2003directsearch} to solve problem~\eqref{eq:pb}. 
Section~\ref{subsec:ddsl} describes a method that uses a Lyapunov function at
every iteration, this variant being close in spirit to decentralized gradient
techniques. Section~\ref{subsec:ddsf} is concerned with an alternative approach 
in which every agent accepts steps solely based on its own function decrease. We point out that both algorithms rely on a different stepsize 
for every agent and every iteration. Our rules for selecting these stepsizes, 
which differ from standard direct-search techniques, will be discussed in detail 
in Section~\ref{sec:cv}.

\subsection{Algorithm based on Lyapunov function decrease}
\label{subsec:ddsl}

Our first algorithm is based on the penalty formulation~\eqref{eq:penpb}, and 
assumes that every agent has access to its own function, neighbor copies of the 
variable, as well as the penalty parameter $\gamma>0$. The objective 
function of~\eqref{eq:penpb} can be rewritten as
\begin{eqnarray*}
	\calL(\bigx;\gamma) &= &\sum_{i=1}^m f_i(x_{i})+\frac{1}{2\gamma}\left( 
	\sum_{i=1}^m \|x_i\|^2-\sum_{i=1}^m \sum_{j=1}^m w_{ij} x_i^\T x_j\right) \\
	&= &\sum_{i=1}^m f_i(x_i)+\frac{1}{2\gamma}\left( \sum_{i=1}^m (1-w_{ii})\|x_i\|^2
	-\sum_{i=1}^m \sum_{j \in \calN_i} w_{ij} x_i^\T x_j\right).
\end{eqnarray*}
When $W$ satisfies Assumption~\ref{as:W} (in particular $\lambda_1(W)=1$), we 
know that the quantity
\[
	\sum_{i=1}^m \left(\|y_i\|^2 - \sum_{j=1}^m w_{ij} y_i^\T y_j \right) 
	= 
	\sum_{i=1}^m \left( (1-w_{ii})\|y_i\|^2 
	- \sum_{j \in \calN_i} w_{ij} y_i^\T y_j  \right)
\]
is nonnegative for any $y_1,\dots,y_m \in \R^n$~\cite{yuan2016convergence}. As a result, under 
Assumption~\ref{as:flow}, $\flow$ is also a lower bound on 
$\calL(\cdot;\gamma)$ for any $\gamma>0$.
Therefore, for a given agent $i$, we consider a specific local Lyapunov function defined as
\begin{equation}
\label{eq:Lifun}
\calL_i\left(x_i;x_{\calN_i},\gamma\right) := 
	f_i(x_i) + \frac{1}{2\gamma}\left[ (1-w_{ii})\|x_i\|^2 - 
	2\sum_{j \in \calN_i} w_{ij} x_i^\T x_j \right],
\end{equation}
where $x_{\calN_i} = \{x_j\}_{j \in \calN_i}$ represents the information transferred 
to agent $i$ from its neighbors. With this definition, we have
\begin{equation*}
	\nabla \calL(\bigx;\gamma) = \left[ 
	\begin{array}{c}
		\nabla_{x_1} \calL_1\left(x_1;x_{\calN_1},\gamma\right) \\
		\vdots \\
		\nabla_{x_m} \calL_m\left(x_m;x_{\calN_m},\gamma\right)
	\end{array}
	\right].
\end{equation*}

During a decentralized optimization process, every 
agent $i$ updates its local copy by performing an (approximate) minimization step of the 
function $\calL_i(\cdot;x_{\calN_i},\gamma)$, then broadcasts its local copy to its 
neighbors and receives their local copies. The agent then updates its function before 
performing another approximate minimization process.

Algorithm~\ref{alg:ddsL} describes a direct-search version of this approach that 
allows every agent to perform one step of direct-search on this local penalty 
function at every iteration. The method follows a standard 
direct-search framework with sufficient decrease that every agent applies in 
parallel to its own Lyapunov function 
$\calL_i\left(x_i;x_{\calN_i},\gamma\right)$. Each agent polls a set of 
directions defined by a positive spanning set $\{D^{(k)}\}$. If there is at least 
one direction for which the local Lyapunov function $\calL_{i}$ exhibits sufficient 
decrease, defined by~\eqref{eq:suffdec}, then a step along that direction 
(scaled by the current stepsize  $\alpha^{(k)}_i$) is taken. Otherwise, this agent's 
local copy is not updated. Communication here is implicit: evaluation
of the conditions for sufficient decrease requires knowledge of the
current estimates of $x^{(k)}_{\calN_i}$, and thus an iteration of 
Algorithm~\ref{alg:ddsL} requires each agent to communicate its local copy to its 
neighbors.

\begin{algorithm}[t]
\caption{Decentralized direct-search based on local Lyapunov decrease (\texttt{DDS-L})}
\label{alg:ddsL}
\begin{algorithmic}
\STATE{\textbf{Inputs:}} Initial points 
$x^{(0)}_{1}=\dots=x^{(0)}_{m}\in \R^n$ and initial stepsizes 
$\alpha_1^{(0)}=\dots=\alpha_m^{(0)}>0$. Consensus 
parameter $\gamma>0$, sequence of positive spanning sets $\{D^{(k)}_{i}\}_{1\le i \le m}$, 
forcing function $\rho: \R^+ \rightarrow \R^+$, mixing matrix $W \in \R^{m \times m}$.

\FOR{each iteration $k=0,1,2,3...$}
\FOR{each agent $i=1,\ldots,m$}
\IF{there exists $d^{(k)}_{i}\in D^{(k)}_i$ such that 
\begin{equation} \label{eq:suffdec}
	\calL_{i}\left(x_{i}^{(k)}+\alpha^{(k)}_i d^{(k)}_{i},x^{(k)}_{\calN_i},\gamma\right) 
	\le \calL_{i}\left(x_{i}^{(k)}, x^{(k)}_{\calN_i},\gamma\right) -\rho(\alpha^{(k)}_i)
 \end{equation}}
\STATE{Set $x^{(k+1)}_{i} = x^{(k)}_{i}+\alpha^{(k)}_i d^{(k)}_{i}$ and declare the iteration  successful for agent $i$.}
\ELSE
\STATE{Set $x^{(k+1)}_{i} = x^{(k)}_{i}$ and declare the iteration unsuccessful 
for agent $i$.}
\ENDIF
\STATE Compute $\alpha^{(k)}_{i+1}$.
\ENDFOR
\ENDFOR
\end{algorithmic}
\end{algorithm}

\subsection{Algorithm based on local function decrease}
\label{subsec:ddsf}

Our second algorithmic proposal is described in Algorithm~\ref{alg:ddsL}. It hews closer to 
standard direct search, in that every agent decides to accept or reject a step based on 
whether a sufficient decrease condition is satisfied for its own local function. This idea 
is also in agreement with the DGD principle~\eqref{eq:dgd}, since the negative gradient of 
the local function $f_i$ but not necessarily for the penalty function.

In addition to differing from Algorithm~\ref{alg:ddsL} in the sufficient decrease acceptance 
condition, Algorithm~\ref{alg:ddsF} also updates the iterate of every agent at each 
iteration through a consensus step involving the mixing matrix $W$. This approach is a 
significant difference from centralized direct-search techniques and is key to guaranteeing 
asymptotic consensus among all agents. As we will show in Section~\ref{subsec:ddsf}, this 
seemingly more natural variant is more challenging to analyze.

\begin{algorithm}[t]
\caption{Decentralized direct-search based on local function decrease (\texttt{DDS-F})}
\label{alg:ddsF}
\begin{algorithmic}
\STATE{\textbf{Inputs:}} Initial points 
$x^{(0)}_{1} =\dots=x^{(0)}_{m}\in \R^n$ and initial stepsizes 
$\alpha_1^{(0)}=\dots=\alpha_m^{(0)}>0$. 
Sequence of positive spanning sets $\{D^{(k)}_i\}$,
forcing function $\rho: \R^+ \rightarrow \R^+$, and mixing matrix 
$W \in \R^{m \times m}$.

\FOR{each iteration $k=0,1,2,3...$}
\FOR{each agent $i=1,\ldots,m$}
\IF{there exists $d^{(k)}_{i}\in D^{(k)}_i$ such that 
\begin{equation} \label{eq:suffdec}
	f_{i}\left(x_{i}^{(k)}+\alpha^{(k)}_i d^{(k)}_{i}\right) 
	\le f_{i}\left(x_{i}^{(k)}\right) -\rho(\alpha^{(k)}_i)
 \end{equation}}
\STATE{Set $x^{(k+1)}_{i} =  \sum\limits_{j\in \mathcal{N}_i \cup \{i\}} w_{ij} x^{(k)}_{j}+\alpha^{(k)}_i d^{(k)}_{i}$ and declare the iteration 
successful for agent $i$.}
\ELSE
\STATE{Set $x^{(k+1)}_{i} = \sum\limits_{j\in \mathcal{N}_i \cup \{i\}} w_{ij} x^{(k)}_{j}$ and declare the iteration unsuccessful 
for agent $i$.}
\ENDIF
\STATE Compute $\alpha^{(k)}_{i+1}$.
\ENDFOR
\ENDFOR
\end{algorithmic}
\end{algorithm}

\section{Convergence results}
\label{sec:cv}

\subsection{Generic assumptions}
\label{subsec:ascv}

This section details the assumptions that are common to our two algorithms.
We will assume, as it is done in classical
directional direct-search~\cite{kolda2003directsearch}, that all positive 
spanning sets considered by the algorithm include bounded directions and have 
cosine measure bounded away from zero.
\begin{assumption}\label{as:pss}
    Consider the sequence $\{D^{(k)}_i\}$ of positive spanning sets used in 
    either Algorithm~\ref{alg:ddsL} or~\ref{alg:ddsF}. There exists 
    $\kappa \in (0,1)$ such that, for each agent $i$ and every 
    $k$, the set $D^{(k)}_i$ is a $\kappa$-descent set, i.e.,
    \begin{equation} \label{eq:centralized:cmDk}
        \cm(D^{(k)}_i) := \min_{v \neq 0_{\real^n}} 
        \max_{d\in {\color{blue}D^{(k)}_i}} \frac{d^\T v}{\|d\|\|v\|} 
        \ge \kappa.
    \end{equation}
    We also assume that all poll directions are normalized, i.e., 
    $\|d\|=1$ for all $d \in D^{(k)}_i$ 
    
\end{assumption}

We note that the choice of using normalized poll directions is made only for the simplicity of exposition, and the analysis easily generalizes to the case where the directions are uniformly bounded in norm, i.e., when there exist 
$0 < \beta_{\min} \le \beta_{\max} < \infty$ such that 
\[
	\forall k,~ \forall d \in D^{(k)}_i, ~~\beta_{\min} \le 
	\|d\| \le \beta_{\max}.
\]
A typical choice of direction set that satisfies Assumption~\ref{as:pss} 
is $D^{(k)}_i=D_{\oplus}=[I\ -I]$ (in that case $\kappa = \tfrac{1}{\sqrt{n}}$).

As in standard convergence analyzes of direct-search schemes, we rely on the 
following requirements for the forcing function. Those guaranty, in particular, 
that the sufficient decrease condition can be satisfied for a sufficiently small 
step size~\cite{KJDzahini_FRinaldi_CWRoyer_DZeffiro_2024}.

\begin{assumption}
\label{as:rho}
	The forcing function $\rho: \R^+ \rightarrow \R^+$ used in either 
    Algorithm~\ref{alg:ddsL} or Algorithm~\ref{alg:ddsF} satisfies the three 
    properties below.
	\begin{enumerate}[(i)]
		\item $\rho$ is nondecreasing,
		\item $\alpha \mapsto \tfrac{\rho(\alpha)}{\alpha}$ is 
		nondecreasing,
		\item $\rho(\alpha)=o(\alpha)$ as $\alpha \rightarrow 0$.
	\end{enumerate}
\end{assumption}

Finally, we state our key assumptions on the step size sequences used by our 
algorithms, which are instrumental in deriving theoretical results for our 
method. Those requirements differ from the classical centralized direct-search setting and aim to cover both adaptive and decreasing stepsize choices. In 
particular, we introduce auxiliary stepsize sequences that control the stepsize 
behavior, akin to predefined stepsize choices in decentralized optimization.

\begin{assumption}
\label{as:alpha}
    For the stepsize sequences $\{\alpha_{i}^{(k)}\}$ used in either 
    Algorithm~\ref{alg:ddsL} or Algorithm~\ref{alg:ddsF}, there exist two 
    sequences $\{\alpha_{\max}^{(k)}\}$ and $\{\alpha_{\min}^{(k)}\}$ such 
    that:
	\begin{enumerate}[(i)]
		\item $\alpha_{\min}^{(k)} \le \alpha_{i}^{(k)} \le \alpha_{\max}^{(k)}$ for all 
		indices $(i,k) \in \{1,\dots,m\} \times \N$;
		\item The sequence $\{\alpha_{\max}^{(k)}\}$ is square summable, i.e.,
		\begin{equation}
		\label{eq:alphaseries}
			\sum_{k \in \N} (\alpha_{\max}^{(k)})^2 < \infty.
		\end{equation}
		\item The sequence $\{\rho(\alpha_{\min}^{(k)})\}$ is not summable, 
		i.e., 
		\begin{equation}
		\label{eq:rhoseries}
			\sum_{k\in \N} \rho(\alpha^{(k)}_{\min})=\infty.
		\end{equation}
	\end{enumerate}
\end{assumption}

Further remarks are in order regarding Assumption~\ref{as:alpha}.
First, square summability of stepsizes as in~\eqref{eq:alphaseries} 
is a standard property that is used in both direct-search~\cite{Gratton_2015} 
and decentralized gradient methods with a unique decreasing
stepsize~\cite{zeng2018nonconvex}.
Secondly, property~\eqref{eq:rhoseries} on the forcing function 
departs from classical choices used in direct-search, such as 
$\rho(\alpha)=\alpha^2$~\cite{kolda2003directsearch}, while relating to 
standard requirements in decentralized gradient techniques.
Thirdly, we assume properties that involve the maximal and minimal stepsizes 
among all agents.

A possible choice for satisfying Assumption~\ref{as:alpha} consists in predefining the sequences $\{\alpha^{(k)}_i\}$ independently of the agents . For instance, for any $k \in \N$, one may set
\begin{equation}
\label{eq:alphabetadec}
	\forall i \in \{1,\dots,m\}, \qquad 
	\alpha_i^{(k)} = \alpha^{(k)}_{\min}=\alpha^{(k)}_{\max} 
	= \frac{\alpha_0}{(1+k)^{\tau_{\alpha}}} 
	\quad \mbox{and} \quad
	\rho(\alpha_i^{(k)}) = \frac{\rho_0}{(1+k)^{\tau_{\rho}}},
\end{equation}
where $\alpha_0>0$, $\rho_0>0$, $0.5 < \tau_{\alpha} < \tau_{\rho} \le 1$. 
It is clear that the resulting sequences satisfy Assumption~\ref{as:alpha}. 
This choice is in line with the standard
practice in decentralized optimization, that favors a priori stepsize 
rules~\cite{nedic2009distributed}.

An alternate stepsize choice, closer to that of direct-search schemes, 
consists in updating $\alpha_i^{(k)}$ in an adaptive fashion for 
every agent. More precisely, if condition~\eqref{eq:suffdec} holds for agent 
$i$ at iteration $k$, then iteration $k$ is successful and one possibly 
increases $\alpha^{(k)}_i$. Otherwise, the iteration is unsuccessful, in which case 
one decreases $\alpha^{(k)}_i$. In addition, the choice of $\rho$ is made 
so as to satisfy Assumption~\ref{as:alpha}.
Overall, for any $i \in \{1,\dots,m\}$ and $k \in \N$, the update formulas are 
given by
\begin{equation}
\label{eq:alphabetaadaptive}
	\alpha^{(k)}_i = \left\{
		\begin{array}{ll}
			\min\left\{\theta^{-1}\alpha^{(k)}_i,\alpha_{\max}^{(k)}\right\} 
			&\mbox{if $k$ is successful for $i$} \\
			\max\left\{\theta\alpha^{(k)}_i,\alpha_{\min}^{(k)}\right\} &\mbox{otherwise,}
		\end{array}
		\right. 
	\quad \mbox{and} \quad 
	\rho(\alpha^{(k)}_i) = c\left(\alpha_i^{(k)}\right)^{1+\tau_{\rho}},
\end{equation}
where $\{\alpha_{\max}^{(k)}\}$ and $\{\alpha_{\min}^{(k)}\}$ are positive sequences 
that converge towards $0$, $\theta \in (0,1)$, $\tau_{\rho} \in (0,1)$, and $c>0$. 
For instance, those sequences can be chosen similarly 
to~\eqref{eq:alphabetadec}, i.e.
\[
    \forall k \in \N, 
    \qquad
    \alpha_{\min}^{(k)} = \frac{c_{\min}}{(1+k)^{\tau_{\alpha}}} 
    \quad \mbox{and} \quad 
    \alpha_{\max}^{(k)} = \frac{c_{\max}}{(1+k)^{\tau_{\alpha}}},
\]
with $0 < c_{\min} < c_{\max}$ and $0.5 < \tau_{\alpha} < \tau_{\rho}$ satisfies 
the desired requirements for convergence. In Section~\ref{sec:num}, we will see that 
the classical strategy in direct search $c_{\min} \rightarrow 0$ and $c_{\max} \rightarrow \infty$ yields good practical performance.

The analysis in the upcoming sections will rely heavily on 
Assumption~\ref{as:alpha} while distinguishing between successful iterations 
(for which at least one agent updates its local copy) and unsuccessful iterations 
(for which all agents leave their local copies unchanged). To this end, we 
let $\calS^{(k)} \subseteq \{1,\dots,m\}$ (resp. 
$\calU^{(k)} \subseteq \{1,\dots,m\}$) denote the set of agents for which 
iteration $k$ is successful (resp. unsuccessful).

\subsection{Convergence analysis of Algorithm~\ref{alg:ddsL} (\texttt{DDS-L})}
\label{subsec:ddsLtheory}

We begin by analyzing Algorithm~\ref{alg:ddsL}.
Since this method relies on a penalty function defined with a constant penalty 
parameter, it cannot be expected to converge to a solution of 
problem~\eqref{eq:conspb}, in the sense that consensus is not guaranteed. 
However, we will show that the method converges towards a stationary point for 
problem~\eqref{eq:penpb}, akin to decentralized gradient 
schemes~\cite{zeng2018nonconvex}.

An iteration of Algorithm~\ref{alg:ddsL} corresponds to applying one step of a 
direct-search algorithm to the function $\calL_i$ for agent $i$. Classical 
analyses of direct-search methods rely on a link between the gradient of the 
objective and the stepsize on unsuccessful iterations. The following lemma 
provides an analogous result for the decentralized setting.

\begin{lemma} 
\label{lem:unsucc}
	Let Assumptions~\ref{as:lip} and \ref{as:pss} hold.
	Suppose that the $k$-th iteration of Algorithm~\ref{alg:ddsL} is unsuccessful 
	for agent $i$. 
	Then, one has
	\begin{equation} 
	\label{eq:unsucc}
		\left\|\nabla_{x_i} \calL_i \left(x^{(k)}_i;x^{(k)}_{\calN_i},\gamma\right)\right\| 
		\le \frac{1}{\kappa}\left(\frac{M_i}{2}\alpha^{(k)}_i + 
		\frac{\rho(\alpha^{(k)}_i)}{\alpha^{(k)}_i} \right),
	\end{equation}
	where $M_i:=L_i+\frac{1-w_{ii}}{\gamma}$.
\end{lemma}

\begin{proof}
	Since the $k$-th iteration is unsuccessful for agent $i$, 
	condition~\eqref{eq:suffdec} does not hold. Therefore, we must have
	\begin{equation}
	\label{eq:nosuffdec}
		\calL_i \left(x^{(k)}_i;x^{(k)}_{\calN_i},\gamma\right) - \rho(\alpha^{(k)}_i)
		< \calL_i \left(x^{(k)}_i + \alpha^{(k)}_i d ; x^{(k)}_{\calN_i},\gamma\right) 
	\end{equation}
	for every $d \in D^{(k)}$. In particular, letting 
	$g^{(k)}_i = \nabla_{x_i} \calL_i \left(x^{(k)}_i;x^{(k)}_{\calN_i},\gamma\right)$ and

	$\bar{d}_i := \arg\max_{d \in D^{(k)}} 
		\frac{d^\T \left[-g^{(k)}_i\right]}{\|d\|\left \|g^{(k)}_i\right\|}$,
	we have by Assumption~\ref{as:pss} that 
	$\bar{d}_i^\T \left[g^{(k)}_i\right] \le -\kappa \|g^{(k)}_i\|$. 
	
	Now, by Assumption~\ref{as:lip}, the function $\calL_i(\cdot;x^{(k)}_{\calN_i},\gamma)$ is 
	continuously differentiable, and its gradient is Lipschitz continuous with Lipschitz 
	constant $M_i^{(k)}:=L_i+\frac{1-w_{ii}}{\gamma}$. As a result,
	\begin{eqnarray*}
		\calL_i \left(x^{(k)}_i + \alpha^{(k)} \bar{d}_i ; x^{(k)}_{\calN_i},\gamma\right) 
		&\le &\calL_i\left(x^{(k)}_i;x^{(k)}_{\calN_i},\gamma\right) + \alpha^{(k)}_i
		\bar{d}_i^\T \left[g^{(k)}_i\right] + \frac{M_i}{2}(\alpha^{(k)}_i)^2 \\
		&\le &\calL_i\left(x^{(k)}_i;x^{(k)}_{\calN_i},\gamma\right) - \alpha^{(k)}_i
		\kappa \|g^{(k)}_i\| + \frac{M_i}{2}(\alpha^{(k)}_i)^2.
	\end{eqnarray*}
	Combining the last inequality with~\eqref{eq:nosuffdec} applied at $\bar{d}$ leads 
	to
	\begin{eqnarray*}
		-\rho(\alpha^{(k)}_i) &\le 
		&\calL_i\left(x^{(k)}_i+\alpha^{(k)}_i \bar{d}_i;x^{(k)}_{\calN_i},\gamma\right) 
		-\calL_i\left(x^{(k)}_i;x^{(k)}_{\calN_i},\gamma\right) \\
		&\le 
		&-\alpha^{(k)}_i \kappa \|g^{(k)}_i\| + \frac{M_i}{2}(\alpha^{(k)}_i)^2.
	\end{eqnarray*}
	
	Re-arranging the terms and replacing $\alpha^{(k)}_i$ and $\beta^{(k)}_i$ by 
	their expressions, we obtain:
	\begin{eqnarray*}
		\|g^{(k)}_i\| &\le &\frac{1}{\kappa}\left[ \frac{M_i}{2}\alpha^{(k)}_i 
		+ \frac{\rho(\alpha^{(k)}_i)}{\alpha^{(k)}_i} \right],
	\end{eqnarray*}
	proving the desired result.
\end{proof}

The contrapositive of Lemma~\ref{lem:unsucc} implies that an iteration 
for which~\eqref{eq:unsucc} does not hold is necessarily successful. In 
centralized direct-search, this property is combined with the fact that 
the stepsize sequence goes to zero to produce convergence and complexity 
guarantees~\cite{vicente2013complexity}. In our case, our assumptions on 
the stepsize and forcing function sequences yield the following result.

\begin{theorem}
\label{th:globalcvfix}
	Let Assumptions~\ref{as:W}, \ref{as:pss}, \ref{as:lip}, \ref{as:flow}, 
and \ref{as:alpha} hold. Suppose further that 
	$\{\alpha_i^{(k)}\}_k \rightarrow 0$ for any $i \in \{1,\dots,m\}$. 
	Then, the sequence of iterates $\left\{ \{x^{(k)}_i\}_{i=1}^m \right\}_k$ generated by 
	Algorithm~\ref{alg:ddsL} satisfies
	\begin{equation}
	\label{eq:globalcvfix}
		\liminf_{k\to\infty} 
        \left\| \nabla \calL(\bigx^{\ks},\gamma) \right\| 
		= 0.
	\end{equation}
\end{theorem}
\begin{proof}
    Since $\alpha_{\max}^{(k)} \rightarrow 0$, we know that for any 
    $\epsilon>0$, there exists $K_{\epsilon}$ such 
	that for every $k \ge K_{\epsilon}$, we have
	\begin{equation}
	\label{eq:ineqalphafixcv}
		\max_{1 \le i \le m} \alpha_i^{(k)} \le \alpha_{\max}^{(k)} < 
		\inf\left\{\ \alpha > 0
		\ \middle|\ 
		\frac{\epsilon}{m} \le 
		\frac{1}{\kappa}\left(\frac{\min_{1 \le i \le m}M_i}{2}\alpha + 
		\frac{\rho(\alpha)}{\alpha} \right)
		\ \right\}.
	\end{equation}
    For the purpose of obtaining a contradiction with~\eqref{eq:globalcvfix}, suppose that there exists  $L\in\mathbb{N}$ such that, $\|\nabla \calL(\bigx^{({\ell})},\gamma)\| > \epsilon$ for any {$\ell \ge L$}. In particular, we have 
    $\|\nabla \calL(\bigx^{(k)},\gamma)\| > \epsilon$ for any $k \ge { \widehat{K}_{\epsilon}:=\max \{K_{\epsilon},L\}}$, 
    which further implies
	\begin{equation}
    \label{eq:contradcv}
		\sum_{i=1}^m \left\|\nabla \mathcal {L}_i(x^{(k)}_i;x^{(k)}_{\calN_i},\gamma)\right\|  
		\ge \epsilon,
	\end{equation}
    Combining~\eqref{eq:ineqalphafixcv} and~\eqref{eq:contradcv}, we find that 
    there must exist an agent $i_k$ such that 
	\[
		\left\|\nabla \mathcal {L}_{i_k}(x^{(k)}_{i_k};x^{(k)}_{\calN_{i_k}},\gamma)\right\| 
		\ge \frac{\epsilon}{m}
	\] 
    and thus iteration $k$ will be successful for agent $i_k$.
	Overall, assuming that~\eqref{eq:globalcvfix} does not hold, there exists 
	{$\widehat{K}_{\epsilon}$} such that every iteration of index $k \ge {\widehat{K}_{\epsilon}}$ is 
	successful for at least one agent. 

	For any $k \ge {\widehat{K}_{\epsilon}}$, a Taylor expansion of $\calL(\cdot,\gamma)$ 
	gives
	\begin{eqnarray*}
		\calL(\bigx^{(k+1)};\gamma) - \calL(\bigx^{k};\gamma) 
		&\le &
		\nabla_{\bigx} \calL(\bigx^{(k)};\gamma)^\T (\bigx^{(k+1)}-\bigx^{(k)}) + 
		\frac{M}{2} \|\bigx^{(k+1)}-\bigx^{(k)}\|^2,
	\end{eqnarray*}
	where $M=\max_{1 \le i \le m} M_i$. Recalling that $\bigx^{(k)}=[x^{(k)}_i]_{i=1}^m$ and 
    
	$\nabla_{\bigx} \calL(\bigx^{(k)};\gamma) = 
	[\nabla_{x_i} \calL(x_i^{(k)};x^{(k)}_{\calN_i},\gamma)]_{i=1}^m$ for every $k$, we 
	obtain:
	\begin{eqnarray}
	\label{eq:ineqfixcv1}
		\calL(\bigx^{(k+1)};\gamma) - \calL(\bigx^{(k)};\gamma) 
		&\le &
		\sum_{i=1}^m \left[ 
		\nabla \calL_i(x^{(k)}_i;x^{(k)}_{\calN_i},\gamma)^\T (x^{(k+1)}_i-x^{(k)}_i) 
		+ \frac{M}{2}\|x^{(k+1)}_i - x^{(k)}_i\|^2 
		\right] \nonumber \\
		&= & 
		\sum_{i \in \calS^{(k)}} \left[ 
		\nabla \calL_i(x^{(k)}_i;x^{(k)}_{\calN_i},\gamma)^\T (x^{(k+1)}_i-x^{(k)}_i) 
		+ \frac{M}{2}\|x^{(k+1)}_i - x^{(k)}_i\|^2 
		\right]. \nonumber \\
        & & 
	\end{eqnarray}
	Using now a Taylor expansion of $\calL_i(x^{(k+1)}_i;x^{(k)}_{\calN_i},\gamma)$ for 
	every $i \in \calS^{(k)}$, we have:
	\[
		\calL_i(x^{(k+1)}_i;x^{(k)}_{\calN_i},\gamma) 
		\ge 
		\calL_i(x^{(k)}_i;x^{(k)}_{\calN_i},\gamma) 
		+ \nabla\calL_i(x^{(k)}_i;x^{(k)}_{\calN_i},\gamma)^\T (x^{(k+1)}_i-x^{(k)}_i)  
		- \frac{M_i}{2} \|x^{(k+1)}_i - x^{(k)}_i\|^2,
	\]
	leading to
	\begin{eqnarray}
	\label{eq:ineqfixcv2}
		\nabla \calL_i(x^{(k)}_i;x^{(k)}_{\calN_i},\gamma)^\T (x^{(k+1)}_i-x^{(k)}_i) 
		&\le 
		&\calL_i(x^{(k+1)}_i;x^{(k)}_{\calN_i},\gamma) 
		- \calL_i(x^{(k)}_i;x^{(k)}_{\calN_i},\gamma)  
		+ \frac{M_i}{2} \|x^{(k+1)}_i - x^{(k)}_i\|^2  \nonumber \\
		&\le 
		&-\rho(\alpha^{(k)}_i) + \frac{M_i}{2} \|x^{(k+1)}_i - x^{(k)}_i\|^2,
	\end{eqnarray}
	since the iteration is successful for agent $i$. 
	Plugging~\eqref{eq:ineqfixcv2} into~\eqref{eq:ineqfixcv1} then gives
	\begin{eqnarray}
	\label{eq:ineqfixcv3}
		\calL(\bigx^{(k+1)};\gamma) - \calL(\bigx^{(k)};\gamma) 
		&\le 
		&\sum_{i \in \calS^{(k)}} \left[ -\rho(\alpha^{(k)}_i) 
		+ \frac{M+M_i}{2}\|x^{(k+1)}_i - x^{(k)}_i\|^2 \right] 
		\nonumber \\
		&= 
		&\sum_{i \in \calS^{(k)}} \left[ -\rho(\alpha^{(k)}_i) 
		+ \frac{M+M_i}{2} (\alpha^{(k)}_i)^2 \right]
		\nonumber\\
		&\le 
		&\sum_{i \in \calS^{(k)}} \left[ -\rho(\alpha^{(k)}_i) + M (\alpha^{(k)}_i)^2 \right]
		\nonumber\\
		&\le 
		&\sum_{i \in \calS^{(k)}} \left[ -\rho(\alpha^{(k)}_{\min}) 
		+ M (\alpha^{(k)}_{\max})^2 \right] 
		\nonumber \\
		\calL(\bigx^{(k+1)};\gamma) - \calL(\bigx^{(k)};\gamma) 
		&\le 
		&-\rho(\alpha^{(k)}_{\min}) 
		+ m M (\alpha^{(k)}_{\max})^2,
	\end{eqnarray}
	where the last inequality comes from the fact that the iteration is 
	successful for at least one agent, and at most for all $m$ agents.

	Now, for any $J > {\widehat{K}_{\epsilon}}$, summing~\eqref{eq:ineqfixcv3}, for any 
	$k \in \{{\widehat{K}_{\epsilon}},\ldots, J-1\}$, gives
	\begin{eqnarray*}
		\sum_{k={\widehat{K}_{\epsilon}}}^{J-1} \rho(\alpha^{(k)}_{\min})
		&\le  
		&\calL(\bigx^{({\widehat{K}_{\epsilon}})};\gamma)-\calL(\bigx^{(J)};\gamma) 
		+ m M \sum_{k={\widehat{K}_{\epsilon}}}^{J-1} (\alpha^{(k)}_{\max})^2 
		\\
		&\le 
		&\calL(\bigx^{({\color{red}\widehat{K}_{\epsilon}})};\gamma)-\flow
		+ m M \sum_{k={\widehat{K}_{\epsilon}}}^{J-1} (\alpha^{(k)}_{\max})^2,
	\end{eqnarray*}
	where the last inequality comes from Assumption~\ref{as:flow}.
	As $J \rightarrow +\infty$, the left-hand side goes to $+\infty$ 
	while the right-hand side is finite by Assumption~\ref{as:alpha}.
	Therefore, we obtain a contradiction, from which we conclude 
	that~\eqref{eq:globalcvfix} holds.
\end{proof}

Theorem~\ref{th:globalcvfix} shows that Algorithm~\ref{alg:ddsL} converges to 
a stationary point of the penalized objective~\eqref{eq:penpb}, but does not 
provide consensus guaranties. As we will show in our experiments, consensus for 
this approach does improve as the penalty parameter $\gamma$ decreases.

\subsection{Convergence analysis of Algorithm~\ref{alg:ddsF} (\texttt{DDS-F})}
\label{subsec:cvddsF}

We now turn to Algorithm~\ref{alg:ddsF}, whose analysis requires a number of 
ingredients from the decentralized optimization literature. On one hand, 
Assumption~\ref{as:alpha} is critical for convergence as it guarantees 
convergence of the stepsize sequence to zero, akin to standard analyses of 
decentralized gradient descent in a nonconvex 
setting~\cite{yuan2016convergence}. 
On the other hand, the properties of the mixing matrix $W$ are instrumental 
for reaching asymptotic consensus, as they yield the following guarantees.

\begin{proposition}\cite[Proposition 1]{nedic2009distributed}
\label{prop:mixmatcons}
Under Assumption~\ref{as:W}, there exists a constant $C_W \ge 0$ such that
\[
\left\|\widehat W^j-\Am \right\|\le C_W \zeta^k
\]
for any $j \in \N$, where $\Am:= \frac{1}{m} \left(\mathbf{1}_m\mathbf{1}_m^\T\right) \otimes \mathbf{I}_n$
and $\zeta$ is the spectral mixing matrix constant defined by
\begin{equation}\label{eq:spect}
\zeta := \max\left(\left\vert \lambda_2(W)\right\vert,
\left\vert \lambda_n(W)\right\vert\right). 
\end{equation}
\end{proposition}

\begin{lemma}\cite[Lemma 8]{zeng2018nonconvex} \label{lem:sumseqzero}
Under Assumptions~\ref{as:W} and~\ref{as:alpha}, there exists $C_{\zeta} \geq 0$ such that
\[
\sum\limits_{j=0}^k \zeta^{k-j} \alpha^{(j)}_{\max}\le C_{\zeta}\alpha^{(k)}_{\max}
\]
for any $k \in \N$.
\end{lemma}

The previous results allow to bound the consensus error at every iteration, by a reasoning 
similar to the derivative-based setting~\cite[Proposition 3]{zeng2018nonconvex}. Since our 
approach does not necessarily update each local copy at every iteration, we provide a full 
proof of that result.
\begin{proposition}\label{prop:cons}
Let Assumption~\ref{as:W} hold. For any $k \in \N$, we have
\[
\left\|\bigx^{(k)}-\Am \bigx^{(k)}\right\|\le 
\sqrt{m} C_W \left(\|\bigx^0\|\zeta^k+C_{\zeta} \alpha^{(k-1)}_{\max}\right) 
\]
where $C_W$ and $C_\zeta$ are the constants defined in Proposition~\ref{prop:mixmatcons} 
and Lemma~\ref{lem:sumseqzero}, respectively.
\end{proposition}
\begin{proof}
By definition of the $k$th iteration, we have 
\[
    \bigx^{\ks} 
    =
    \widehat W^k \bigx^{(0)} 
    + 
    \sum_{j=0}^{k-1} \widehat W^{k-1-j} \left(\alpha^{(j)} \otimes \mathbf{1}_n\right) 
    \odot \bigd^{(j)},
\]
where $\bigd^{(j)}\in\mathbb{R}^{mn}$ concatenates the directions used by each agent at 
iteration $j$, i.e. $d_i^{(j)}$ if $i \in \mathcal{S}^{(j)}$ and $0_{\R^n}$ otherwise. 
As a result,
\begin{eqnarray*}
    \bigx^{\ks} - \Am\,\bigx^{(k)} 
    &=
    & \left(\widehat W^k - \Am \widehat W^k\right) \bigx^{(0)} + \sum_{j=0}^{k-1} \left(\widehat W^{k-1-j} - \Am
    \widehat W^{k-1-j}\right)\left(\alpha^{(j)} \otimes \mathbf{1}_n\right) 
    \odot \bigd^{(j)}.
\end{eqnarray*}
Meanwhile, by Assumption~\ref{as:W}, we have $\Am \widehat W^{k-1-j} \bigv = \Am \bigv$
for any $\bigv \in \R^{nm}$ and any $j \in \{0,\dots,k-1\}$. Using this property 
together with Cauchy-Schwarz inequality, we obtain
\begin{eqnarray*}
    \left\| \bigx^{\ks}-\Am \bigx^{\ks} \right\| 
    &\le 
    &\left\| \widehat W^k - \Am \right\| \|\bigx^{(0)}\| 
    + \sum_{j=0}^{k-1} \|\widehat W^{k-1-j} - \Am\| 
    \left\| \left(\alpha^{(j)} \otimes \mathbf{1}_n\right) 
    \odot \bigd^{(j)} \right\|.
\end{eqnarray*}
For any $j \in \{0,\dots,k-1\}$, we have
\[
    \left\| \left(\alpha^{(j)} \otimes \mathbf{1}_n\right) \odot \bigd^{(j)} \right\| 
    = \sqrt{\sum_{i \in \mathcal{S}^{(j)}} \left(\alpha_i^{(j)} \|d_i^{(j)}\| \right)^2}
    \le \sqrt{m}\alpha_{\max}^{(j)},
\]
where we used that $\|\bigd^{(j)}\|^2 = \sum_{i \in \calS^{(j)}} \|d_i^{(j)}\|^2$ and 
$\|d_i^{(j)}\|=1$. Thus,
\begin{eqnarray*}
    \left\| \bigx^{\ks}-\Am \bigx^{\ks} \right\|  
    &\le 
    &\left\| \widehat W^k - A_m \right\| \|\bigx^{(0)}\| 
    +\sum_{j=0}^{k-1} \|\widehat W^{k-1-j} - \Am\|  \sqrt{m}\alpha_{\max}^{(j)} \\
    &\le 
    &C_W \zeta^k \|\bigx^{(0)}\| 
    +\sum_{j=0}^{k-1} \sqrt{m} C_W \zeta^{k-1-j} \alpha_{\max}^{(j)} \\
    &\le 
    &C_W \|\bigx^{(0)}\| \zeta^k 
    + C_W \sqrt{m} \alpha_{\max}^{(k-1)}
\end{eqnarray*}
where the second inequality follows from Proposition~\ref{prop:mixmatcons} and the 
third one follows from Lemma~\ref{lem:sumseqzero}. Rearranging the constants yields 
the desired conclusion.
\end{proof}

Provided the stepsize is chosen to satisfy our assumptions, one can establish that 
the iterates of Algorithm~\ref{alg:ddsF} reach asymptotic consensus.

\begin{theorem}
\label{th:consensusDDSF}
Under Assumptions~\ref{as:W} and~\ref{as:alpha}, the iterates of Algorithm~\ref{alg:ddsF} 
satisfy
\[
    \lim_{k \rightarrow \infty} \left\|\bigx^{(k)}-\widehat W \bigx^{(k)}\right\| = 0.
\]
\end{theorem}
\begin{proof}
For any $k \in \N$, we have
\[
    \|\bigx^{\ks} - \widehat W \bigx^{\ks} \| 
    \le 
    \|\bigx^{\ks} - \Am \bigx^{\ks}\| + \|\Am \bigx^{\ks}-\widehat W \bigx^{\ks}\|.
\]
Noticing that 
\[
    \Am \widehat W \bigx^{\ks} 
    = \Am \bigx^{\ks}
    = \widehat W \Am \bigx^{\ks},
\]
we obtain
\begin{eqnarray}
\label{eq:consensusDDFproof}
    \|\bigx^{\ks} - \widehat W \bigx^{\ks} \| 
    &=
    &\|\bigx^{\ks} - \Am \bigx^{\ks}\| 
    + \|\widehat W \Am \bigx^{\ks}-\widehat W \bigx^{\ks}\| \nonumber \\
    &\le
    &(1+ \|\widehat W\|) \|\bigx^{\ks} - \Am \bigx^{\ks}\| \nonumber \\
    &\le 
    &\sqrt{m} C_W (1+ \|\widehat W\|)\left(
    \|\bigx^0\|\zeta^k+C_{\zeta} \alpha^{(k-1)}_{\max}\right),
\end{eqnarray}
where the last inequality is from Proposition~\ref{prop:cons}. By Assumption~\ref{as:W}, 
$\zeta \in (0,1)$ and thus $\lim_{k \rightarrow \infty} \zeta^k = 0$. By 
Assumption~\ref{as:alpha}, $\sum_k (\alpha_{\max}^{(k)})^2 < \infty$ and therefore 
\[
    \lim_{k \rightarrow \infty} \alpha_{\max}^{(k)} 
    =
    \lim_{k \rightarrow \infty} \alpha_{\max}^{(k-1)} = 0. 
\]
Combining these observations with~\eqref{eq:consensusDDFproof} yields the 
desired conclusion.
\end{proof}

The result of Theorem~\ref{th:consensusDDSF} is somewhat expected, since 
Algorithm~\ref{alg:ddsF} performs a consensus step at every iteration for all 
agents, regardless of the successful nature of that iteration. Note however 
that Assumption~\ref{as:alpha} (and in particular the fact that 
$\{\alpha_{\max}^{(k)}\}$ converges to zero) is instrumental to such a result.

Unlike DGD techniques, however, we cannot establish convergence to a first-order 
stationary point for this framework, in part because the decrease condition may 
not be in line with the consensus step, while relying on directions that are not 
directly related to the true gradients. We illustrate this issue using a 
two-dimensional example with two agents. Suppose that \( n = m = 2 \) and that  
\[
    f_1(x) = (x_1 - 1)^2, \quad f_2(x) = x_2^2 
    \quad \mbox{and} \quad
    W = \frac{1}{2} \begin{bmatrix} 1 & 1 \\ 1 & 1 \end{bmatrix}.
\]
Suppose further that we apply Algorithm~\ref{alg:ddsF} with  
\( x_1^{(0)} = x_2^{(0)} = \begin{bmatrix} 0 & 1 \end{bmatrix}^T \) and  
\( D_k = D = \{d, -d, \dots\} \) with \( d = \begin{bmatrix} 1 & 1 \end{bmatrix}^T \). 
Then, any decreasing stepsize sequence \( \alpha^{(k)} \) such that 
the iteration is always successful for agent 1 using direction \( d \) and for 
agent 2 using direction \( -d \) would lead to consensus without optimality. 
Indeed, in such a case, the average iterate would always equal
\( \begin{bmatrix} 0 & 1 \end{bmatrix}^T \), while the iterates would have 
the form  
\[
    x_1^{(k)} = \begin{bmatrix} \alpha_k \\ 1 + \alpha_k \end{bmatrix}, \quad  
    x_2^{(k)} = \begin{bmatrix} -\alpha_k \\ 1 - \alpha_k \end{bmatrix},
\]
implying that Algorithm~\ref{alg:ddsF} never converges to a minimizer. This 
example casts doubts on the convergence guaranties of Algorithm~\ref{alg:ddsF}, 
yet we have found that method to perform quite well in our experiments, 
described in the next section.

\section{Numerical Results}
\label{sec:num}

In this section, we evaluate the performance of our proposed direct-search 
algorithms compared to zeroth-order variants of decentralized gradient descent. 
We first compare our algorithms on a toy problem from the decentralized 
literature~\cite{hajinezhad2017zeroth}. We then adapt the Mor\'e-Wild 
test set~\cite{JJMore_SMWild_2009}, a standard benchmark in derivative-free 
optimization, to the decentralized setting.

\subsection{Implementation details}
\label{ssec:implem}

Our experiments are conducted in MATLAB version R2024a. 
Algorithms~\ref{alg:ddsL} (\texttt{DDS-L)} and~\ref{alg:ddsF} 
(\texttt{DDS-F}) were implemented using the same positive spanning set 
across all iterations, namely, for every agent  $i$, ${D^{(k)}_i=D_{\oplus}}=[B_{\oplus}, -B_{\oplus}]$ where $B_{\oplus}$ 
corresponds to the canonical basis of $\mathbb{R}^n$. For each method, we 
considered two variants based on the following stepsize updating rules:
\begin{itemize}
    \item \texttt{Vanishing}:  $\alpha^{(k)}=\alpha^{(k)}_{\max}=\alpha^{(k)}_{\min}=\frac{\alpha^0}{(1+k)^{0.6}}$, $c=10^{-8}$, and $\tau_{\rho}=0.8$;
    \item \texttt{Adaptive}: $\alpha^{(k)}_{\max}=+\infty$, 
    { $\alpha^{(k)}_{\min}=0$}, $\theta=0.5$, $c=10^{-8}$, and $\tau_{\rho}=0.8$.
\end{itemize}

We compared \texttt{DDS-L} and \texttt{DDS-F} with two zeroth-order 
decentralized gradient schemes based on the iteration
\begin{equation}
\label{eq:zodgd}
    x_{i}^{\k1s} = \sum_{j \in \calN_i} w_{ij} x_{j}^\ks - \alpha^{(k)} \tilde{g}_i^\ks 
    \qquad 
    \forall k \in \N,\ \forall i \in \{1,\dots,m\},
\end{equation}
where $\tilde{g}_i^{(k)}$ is a gradient approximation built from function values. The variant 
\texttt{ZO-DGD(FD)} is built from centered finite differences based on evaluating $f_i$ 
along all directions in $D$ with a finite difference parameter of $10^{-7}$. The variant 
\texttt{ZO-DGD(LM)} fits a linear model to previously available values, in the spirit of 
model-based derivative-free optimization~\cite{conn2009introduction}.

We run all solvers using $\alpha^{0}=\|x_0\|+1$ as an initial stepsize, where $x_0$ is the same starting point for all solvers. The underlying network for all problems is a graph 
$\mathcal{G}=(\mathcal{V},\mathcal{E})$ of $m$ agents, where agents are connected with 
a probability $p_c=0.5$.

At each iteration $k$, we record the following metrics:
\begin{itemize}
    \item $\sum_{i=1}^m f_i (x^{(k)}_{(i)})$, where $x^{(k)}_{(i)}$ is the iterate associated with agent $i$ at iteration $k$,
    \item $\sum_{i=1}^m f_i (\bar{x}^{(k)})$, where $\bar{x}^{(k)}:= \frac{1}{m}\sum_{i=1}^m x^{(k)}_{(i)}$ is the average of the iterates associated with all agents $m$ at iteration $k$,
    \item $\sum_{i=1}^m \|x^{(k)}_{(i)} -\bar{x}^{(k)}\|$, representing the consensus of all $m$ agents at iteration $k$.
\end{itemize}
The first two metrics are associated with optimality, while the last metric measures
the agreement between the local iterates.

\subsection{A separable problem}
\label{ssec:toypb}

We consider the following optimization problem from Hajinezhad et al~\cite{hajinezhad2017zeroth}:
\begin{equation}
\label{eq:toypb}
\min_{x \in \real^n} \sum_{i=1}^n f_i(x_i)= \frac{a_i}{1 + \exp(-x_i)} + b_i \log(1+ x_i^2) ~~~~~~~~~ x_i = x_j \quad \forall (i,j) \in \mathcal{V}.
\end{equation}
Problem~\eqref{eq:toypb} fits our general framework~\eqref{eq:pb} with $m=n$. Moreover, 
note that the problem is \emph{separable}, in that every function $f_i$ depends solely 
on the $i$th optimization variable. We consider three instances of problem~\eqref{eq:toypb} 
corresponding to $n \in \{5,10,15\}$. In each instance, the parameters 
$\{a_i,b_i\}_{i=1,\dots,m}$ were generated following an i.i.d. standard Gaussian 
distribution. All methods were run with a maximum budget of $100n$ evaluations, using the 
vector of all ones in $\R^n$ as a starting point.

\begin{figure}[!ht]
\centering
\subfigure[$\gamma=1$ \label{cv:fig:toy:pb:5:1}]{
\includegraphics[ scale=0.275]{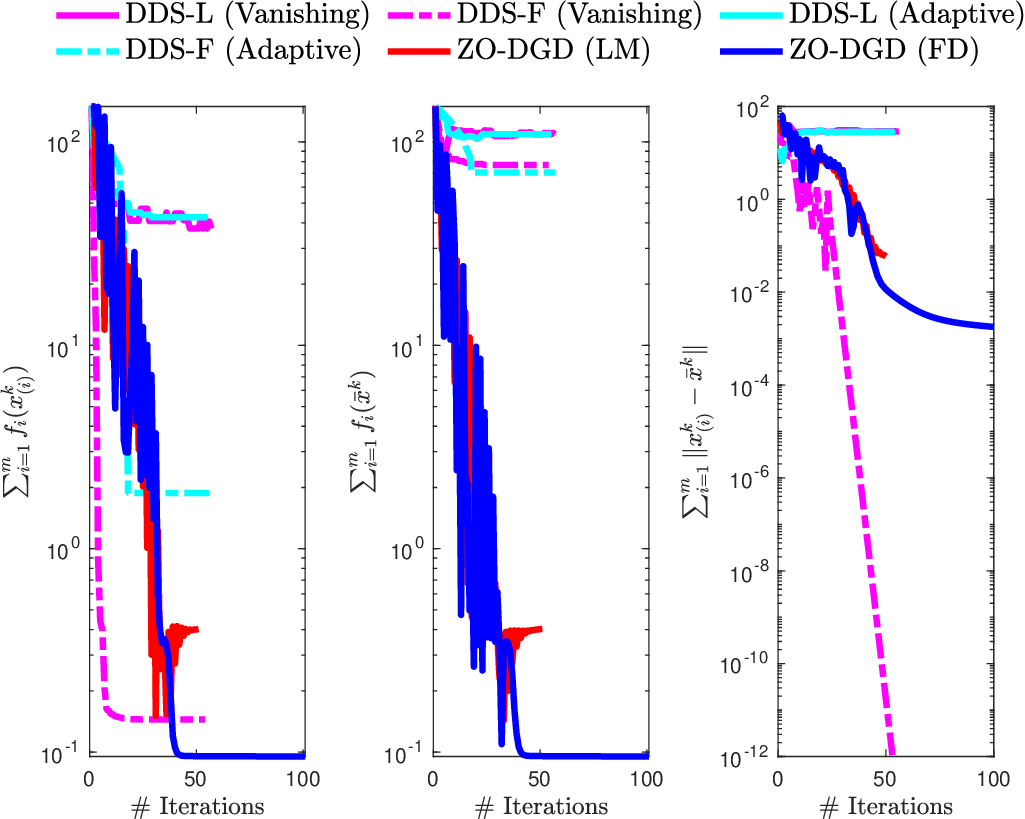}
}
\subfigure[$\gamma=10$ \label{cv:fig:toy:pb:5:1}]{
\includegraphics[scale=0.275]{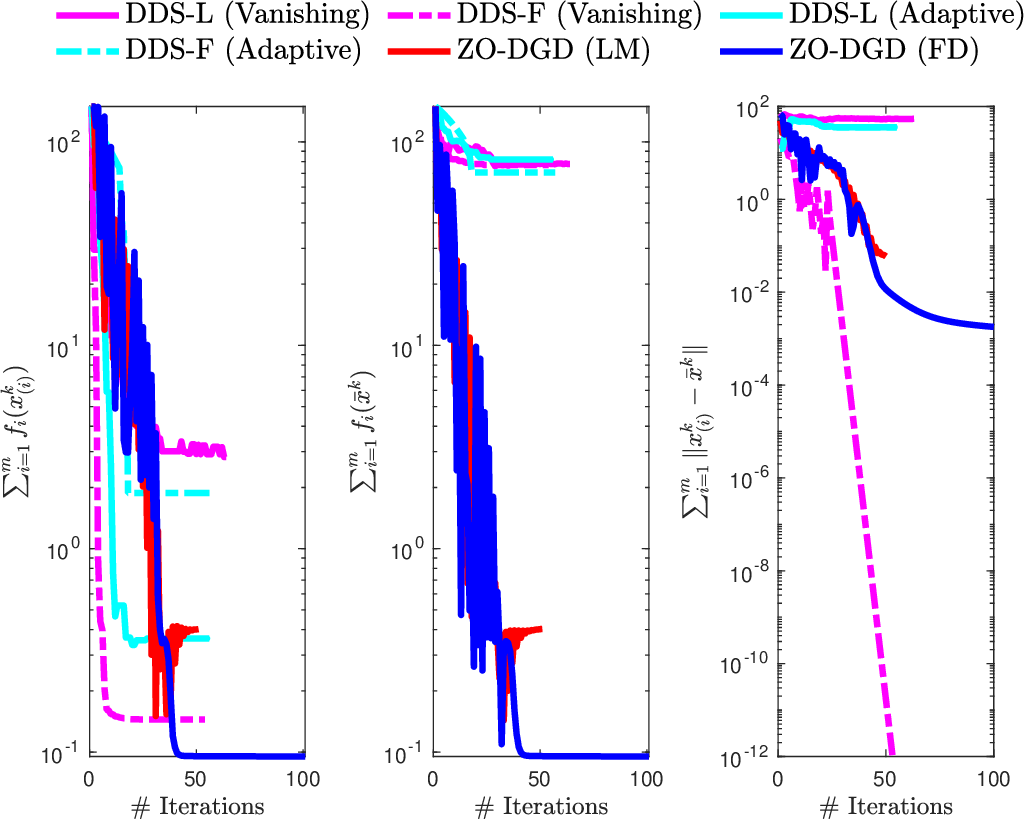}
}
\subfigure[$\gamma=100$ \label{cv:fig:toy:pb:5:10}]{
\includegraphics[scale=0.275]{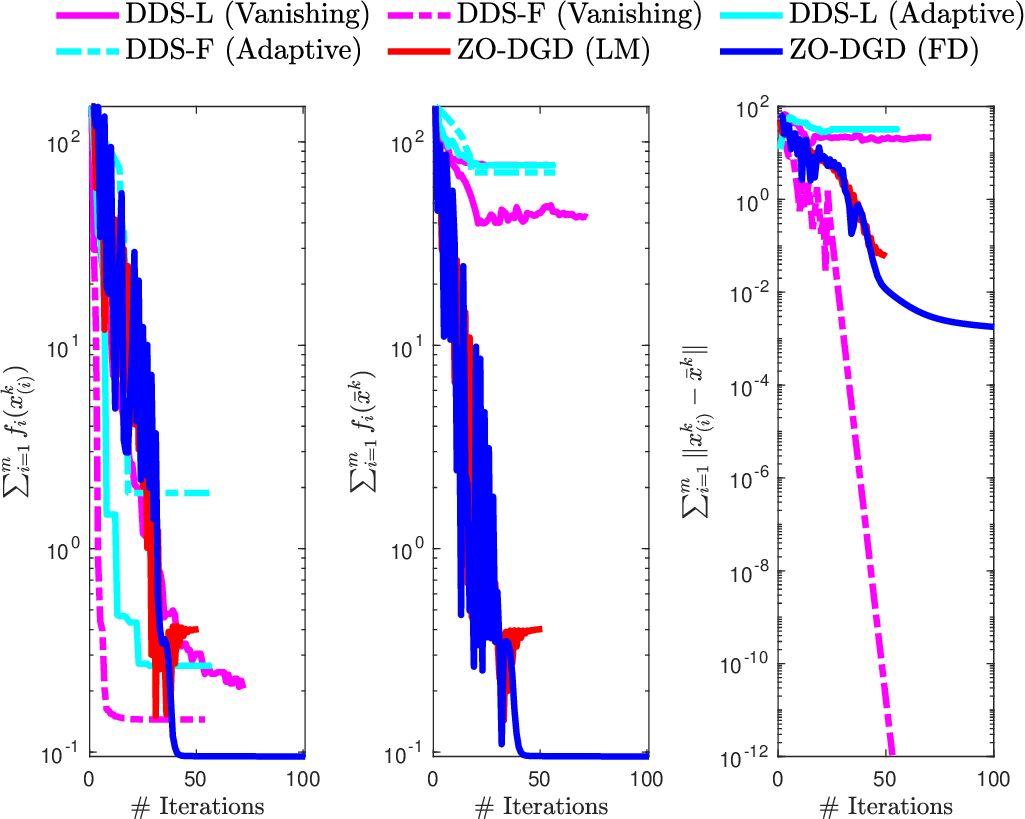}
}
\caption{Convergence plots for problem~\eqref{eq:toypb} in dimension $n=5$.}
\label{fig:toypb:2:n:5}
\end{figure}

Figure~\ref{fig:toypb:2:n:5} focuses on the case $n=5$. Our goal is both to compare our 
algorithms with zeroth-order schemes and to investigate the impact of $\gamma$ on the 
performance of \texttt{DDS-L} (Algorithm~\ref{alg:ddsF}). First, note that the best 
variant in terms of objective value at the iterates and the averaged iterate is the 
finite-difference variant \texttt{ZO-DGD (FD)}. However, we note that \texttt{DDS-F} with 
vanishing stepsizes converges more quickly in terms of objective value, even though it 
plateaus at a higher value overall. Besides, the \texttt{DDS-F (Vanishing)} variant 
outperforms the other methods in terms of consensus, which is a key aspect in decentralized 
algorithms. We note that increasing $\gamma$ improves consensus for the \texttt{DDS-L} 
variants, although those variants are outperformed by the zeroth-order schemes as well 
as the \texttt{DDS-F} ones.

\begin{figure}[!ht]
\centering
\subfigure[$n=m=5$ \label{cv:fig:toy:pb:2:n:5}]{
\includegraphics[scale=0.275]{./Figures/f23_connectproba_0.5_gamma_1.eps}
}
\subfigure[$n=m=10$ \label{cv:fig:toy:pb:2:n:10}]{
\includegraphics[scale=0.275]{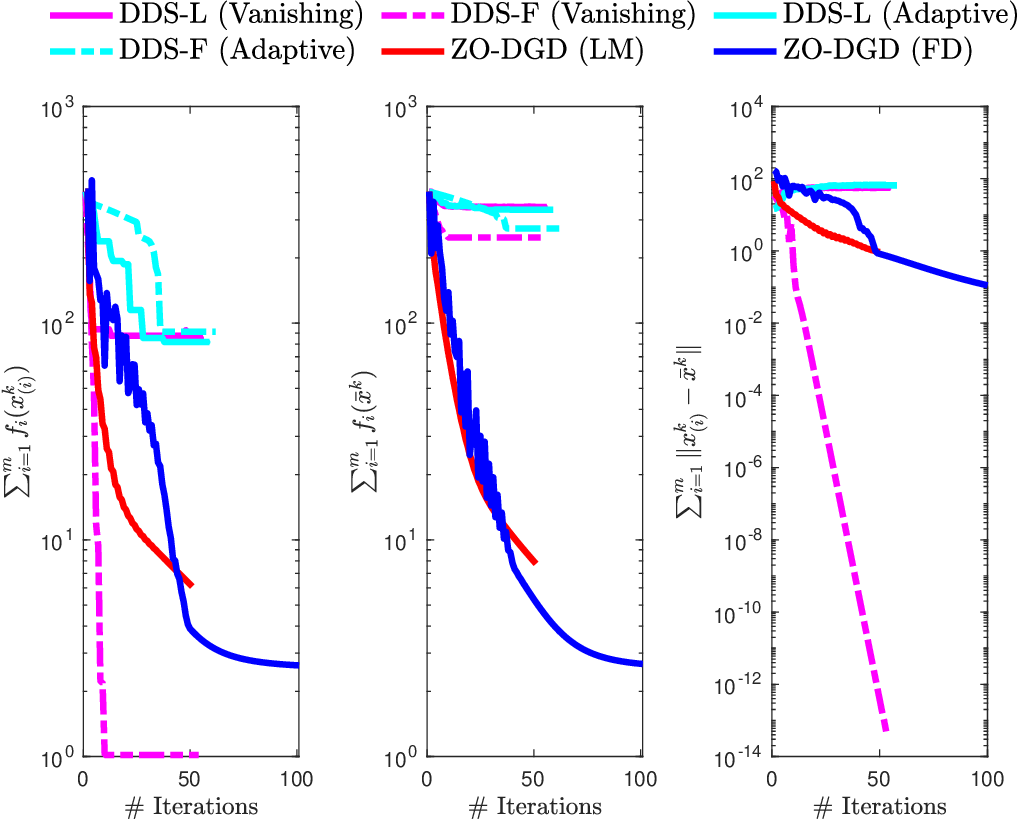}
}
\subfigure[$n=m=15$ \label{cv:fig:toy:pb:2:nb:15}]{
\includegraphics[scale=0.275]{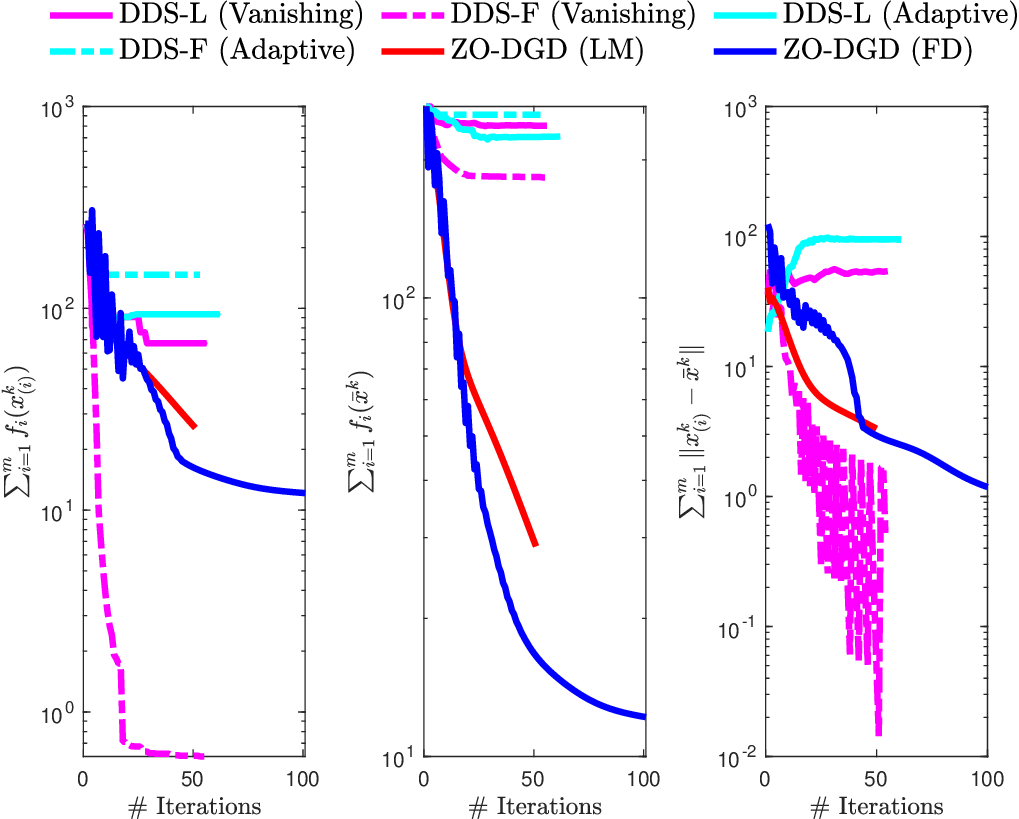}
}

\caption{Convergence plots for problem~\eqref{eq:toypb} with $\gamma=1$.}
\label{fig:toypb:2}
\end{figure}

Figure~\ref{fig:toypb:2} illustrates the variability in performance as the dimension of the 
problem changes. We observe again that \texttt{DDS-F (Vanishing)} outperforms the other 
variants in terms of objective and consensus values for most values, while the finite-difference
 zeroth-order scheme yields the lowest average objective value. Overall, these 
results suggest that classical direct-search approaches can outperform zeroth-order schemes 
in a decentralized setting, akin to the centralized case.

\subsection{Decentralized Mor\'e-Wild test set}
\label{ssec:morewild}

We now compare our algorithms on a test set built from that of Mor\'e and 
Wild~\cite{JJMore_SMWild_2009}. This test set comprises 22 nonlinear smooth vector 
functions of the form $F:\R^n \rightarrow \R^m$ where $2 \le n \le 30$ and $2 \le m\le 65$. 
In our experiments, we consider that the components of $F$ are aggregated as a sum of 
squares, yielding an objective of the form $\sum_{i=1}^m F_i(x)^2$. The local function of 
agent $i$ is thus the function $f_i : x \mapsto F_i(x)^2$. We run our algorithms using
the default starting points of the test set~\cite{JJMore_SMWild_2009} as well as a 
maximum budget of either $400 n m$  local function evaluations (i.e., a budget of $400 n$ evaluations per agent ) or $500$ iterations.

The computational analysis is carried out using well-known tools from the literature; that is, performance and data  profiles~(see~\cite{EDDolan_JJMore_2002,JJMore_SMWild_2009} for further details). 
We briefly recall their definitions.	Given a set~$\mathcal{S}$ of algorithms and a set~$\mathcal{P}$ of problems, for $s\in \mathcal{S}$ and $p \in \mathcal{P}$, let $t_{p,s}$ be the number of function evaluations required by algorithm $s$ on problem~$p$ to satisfy the condition 
     \begin{equation*}
     	\texttt{opt}(\bigx^{(k)}) \leq  \texttt{opt}_{\mbox{low}} + \alpha\left( \texttt{opt}(\bigx^{(k)}) -  \texttt{opt}_{\mbox{low}}\right)\, , 
     \end{equation*}     
     where {$\alpha \in (0, 1)$} , $\texttt{opt}(\bigx^{(k)})$ is the optimality metric (i.e.,  $\sum_{i=1}^m f_i (x^{(k)}_{(i)})$, $\sum_{i=1}^m f_i (\bar{x}^{(k)})$, or $\sum_{i=1}^m \|x^{(k)}_{(i)} -\bar{x}^{(k)}\|$), and $ \texttt{opt}_{\mbox{low}}$ is the best optimality metric value achieved by any solver on problem~$p$. Then, the performance and data profiles of solver $s$ are defined by
     \begin{eqnarray*}
     	\rho_s(\gamma) & := & \frac{1}{|\mathcal{P}|}\left|\left\{p\in \mathcal{P}: \frac{t_{p,s}}{\min\{t_{p,s'}:s'\in \mathcal{S}\}}\leq\gamma\right\}\right|,\\
     	d_s(\kappa) & := & \frac{1}{|\mathcal{P}|}\left|\left\{p\in \mathcal{P}: t_{p,s}\leq\kappa(n_p+1)\right\}\right|\, ,
     \end{eqnarray*}
     where $n_p$ is the dimension of problem $p$. In our tests, for both data and performance profiles, we used two tolerance choices for $\alpha$: $10^{-3}$ and $10^{-6}$.

\begin{figure}[!ht]
\centering
\subfigure[$\alpha=10^{-3}$\label{pp:01:1e-3}]{
\includegraphics[scale=0.3]{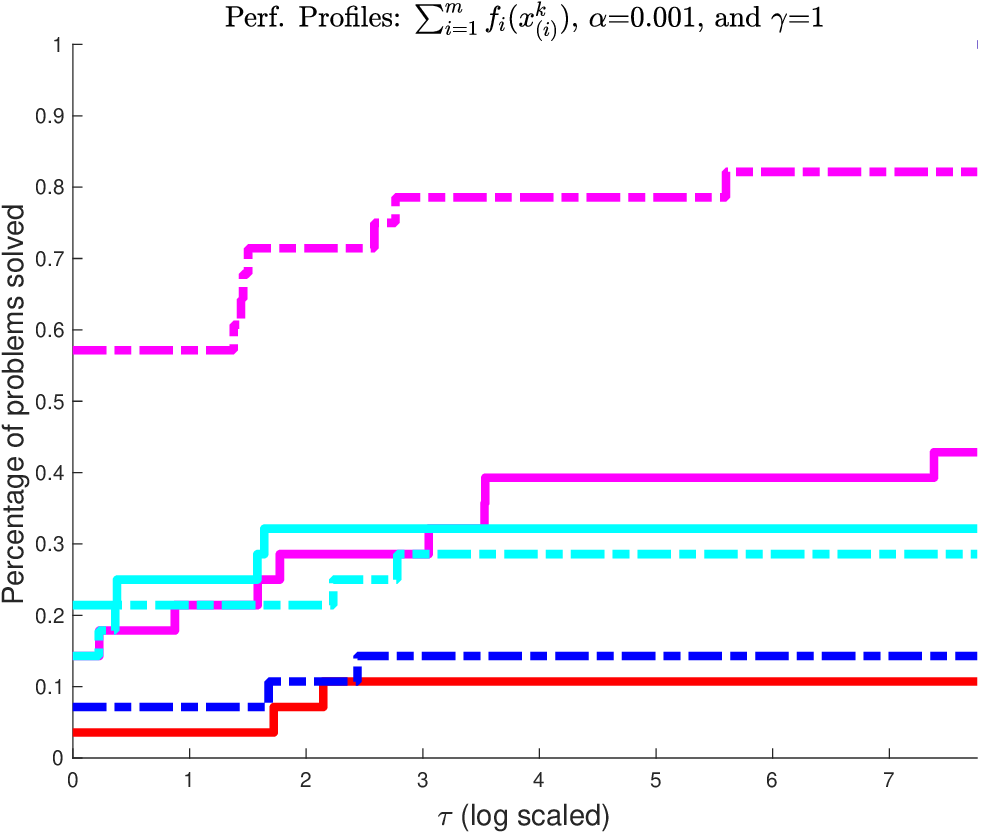}
\includegraphics[scale=0.3]{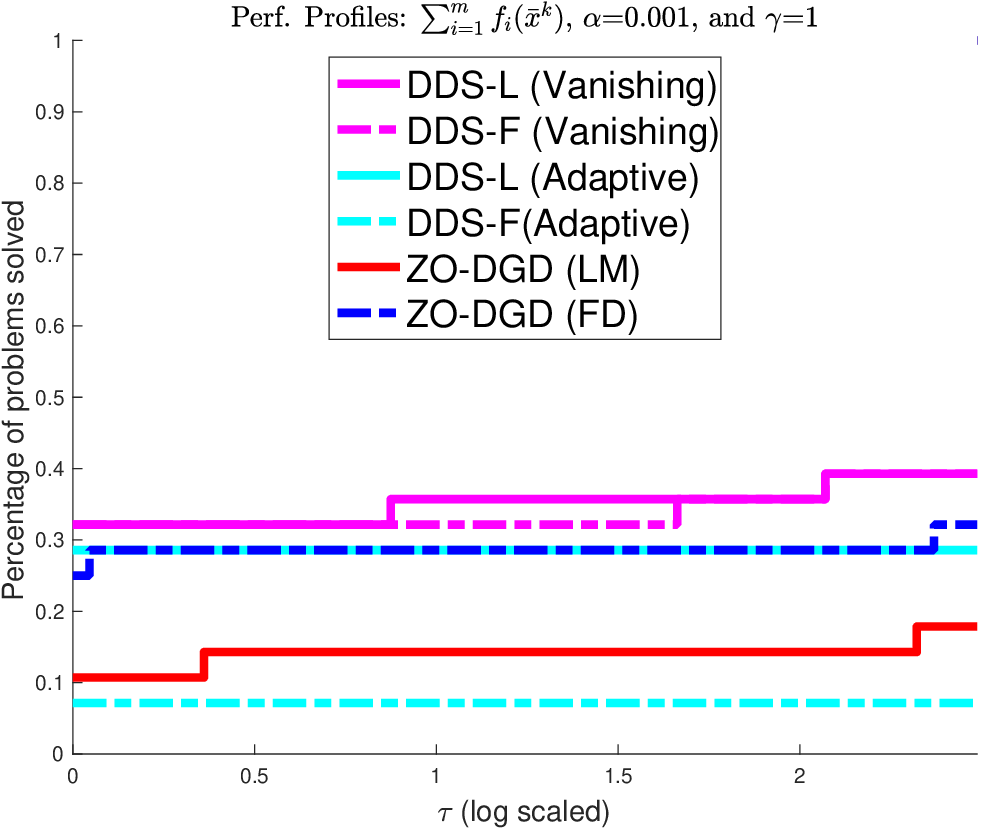}
\includegraphics[scale=0.3]{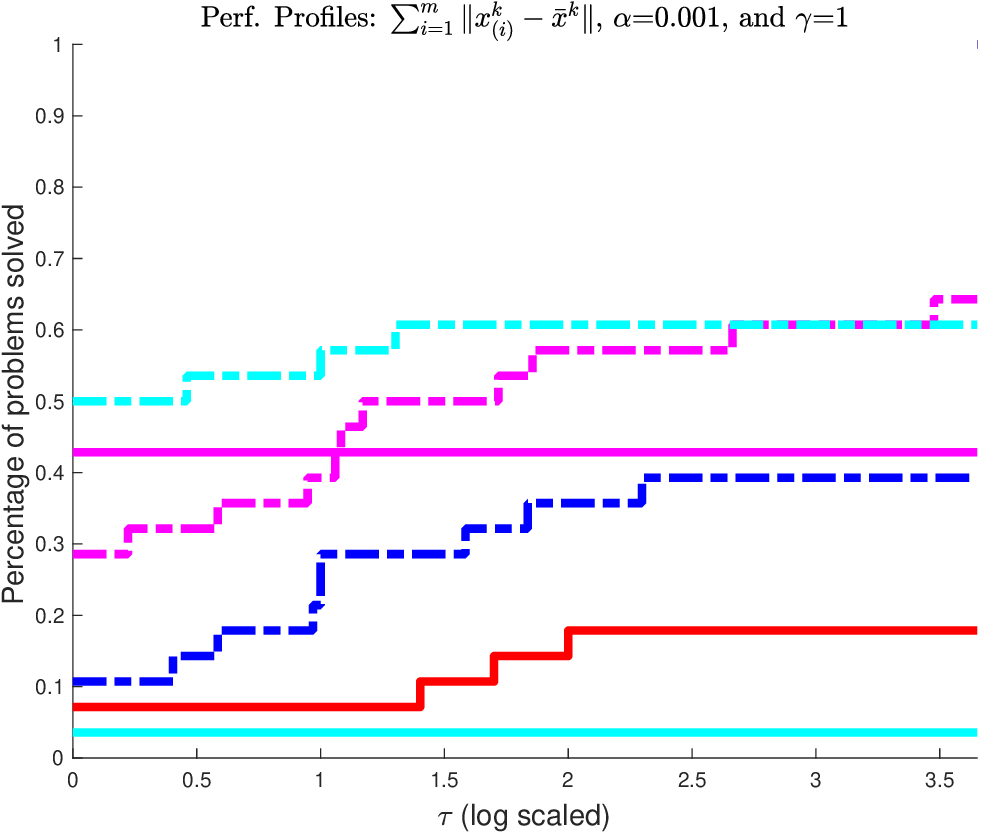}
}
\subfigure[$\alpha=10^{-6}$ \label{pp:01:1e-06}]{
\includegraphics[scale=0.3]{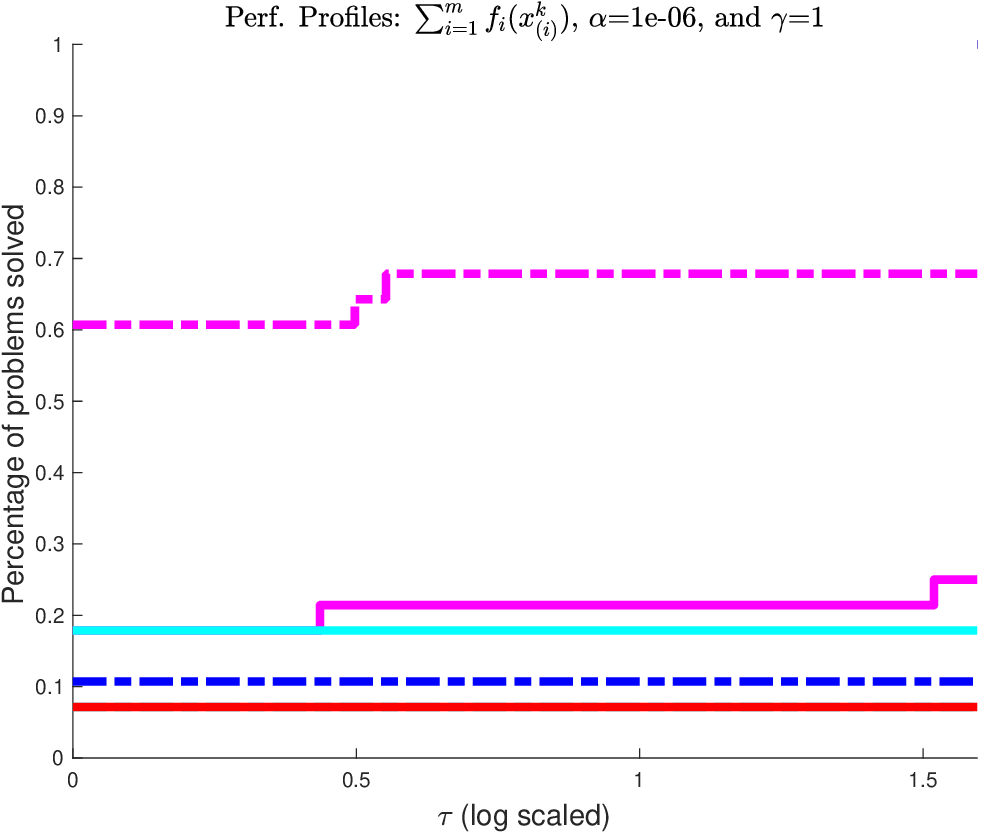}
\includegraphics[scale=0.3]{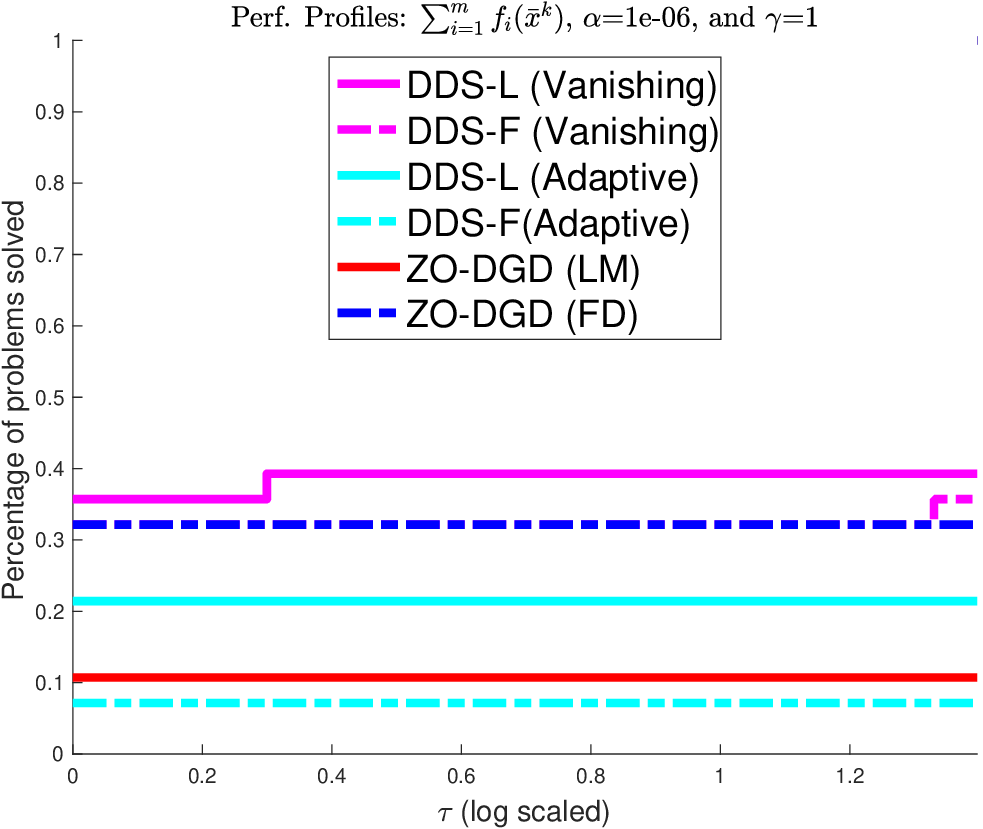}
\includegraphics[scale=0.3]{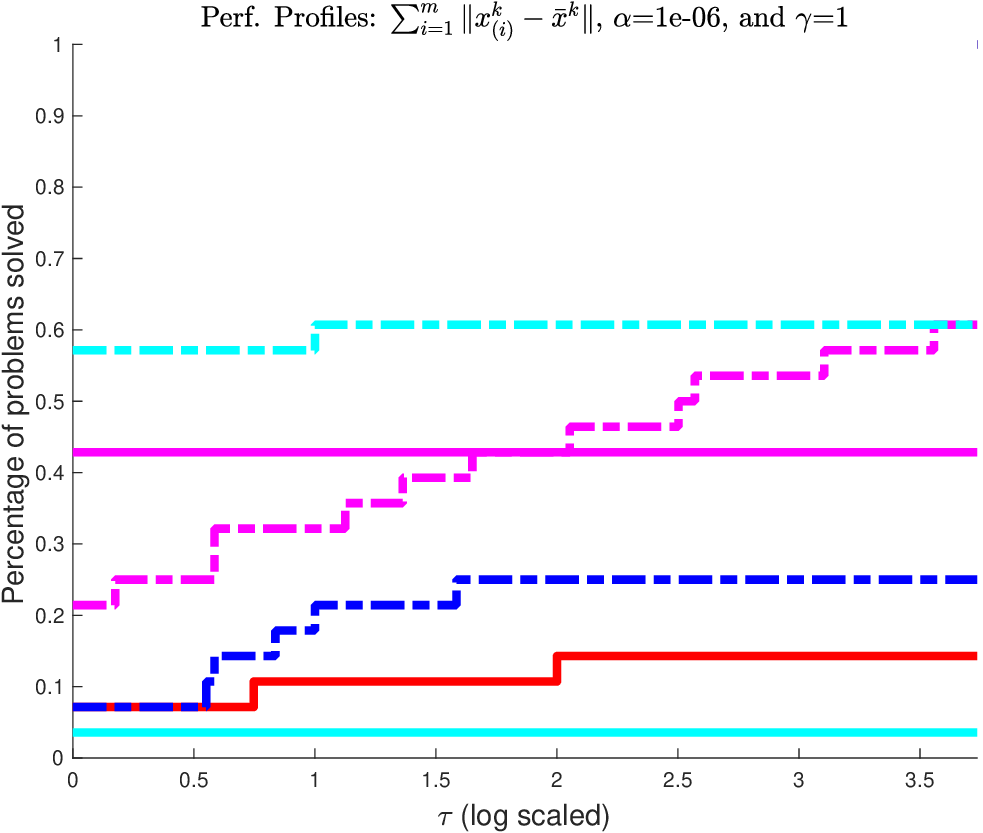}
}
\caption{Performance profiles using three different optimality metrics.}
\label{pp:fig:1e-3}
\end{figure}

Figure~\ref{pp:fig:1e-3} shows performance profiles~\cite{EDDolan_JJMore_2002} 
comparing our various algorithms. We observe that the \texttt{DDS} variants 
mostly outperform the \texttt{ZO-DGD} methods in terms of function values 
and consensus, with \texttt{DDS-F (Vanishing)} standing out as the best 
variant overall. We note that the discrepancy between direct-search and 
zeroth-order methods is less pronounced in terms of function values at the 
average iterates, which is on par with the observations of 
Section~\ref{ssec:toypb}.

\begin{figure}[!ht]
\centering
\subfigure[$\alpha=10^{-3}$ \label{dp:01:1e-3}]{
\includegraphics[scale=0.3]{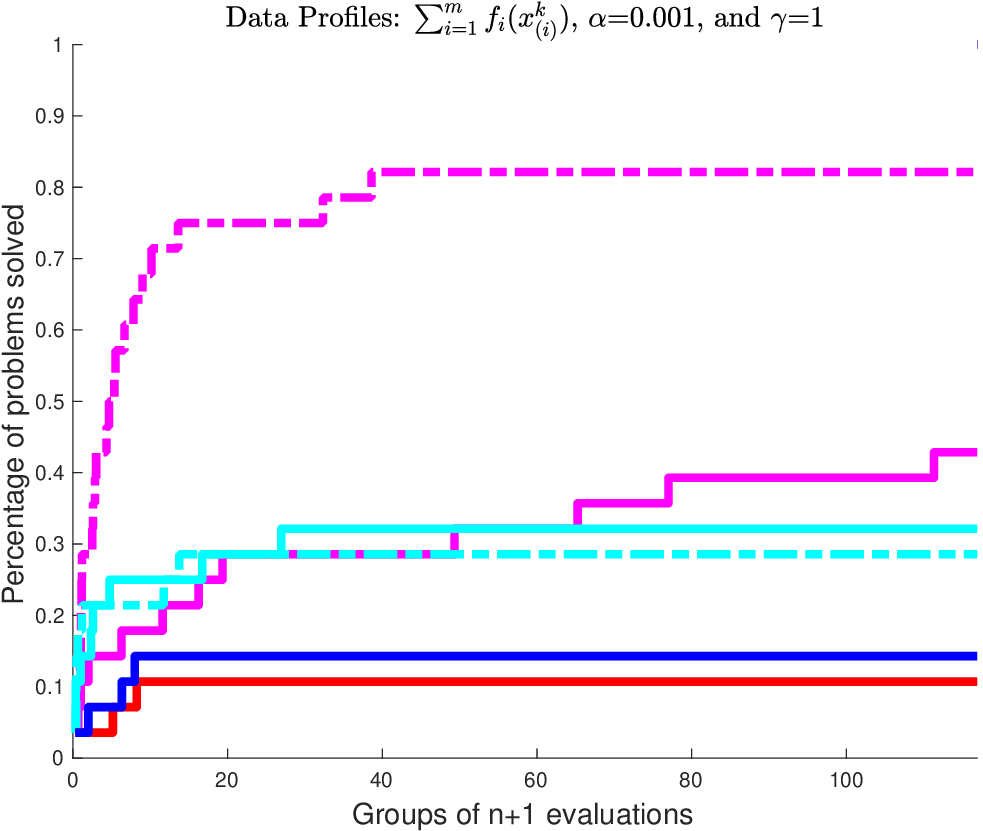}
\includegraphics[scale=0.3]{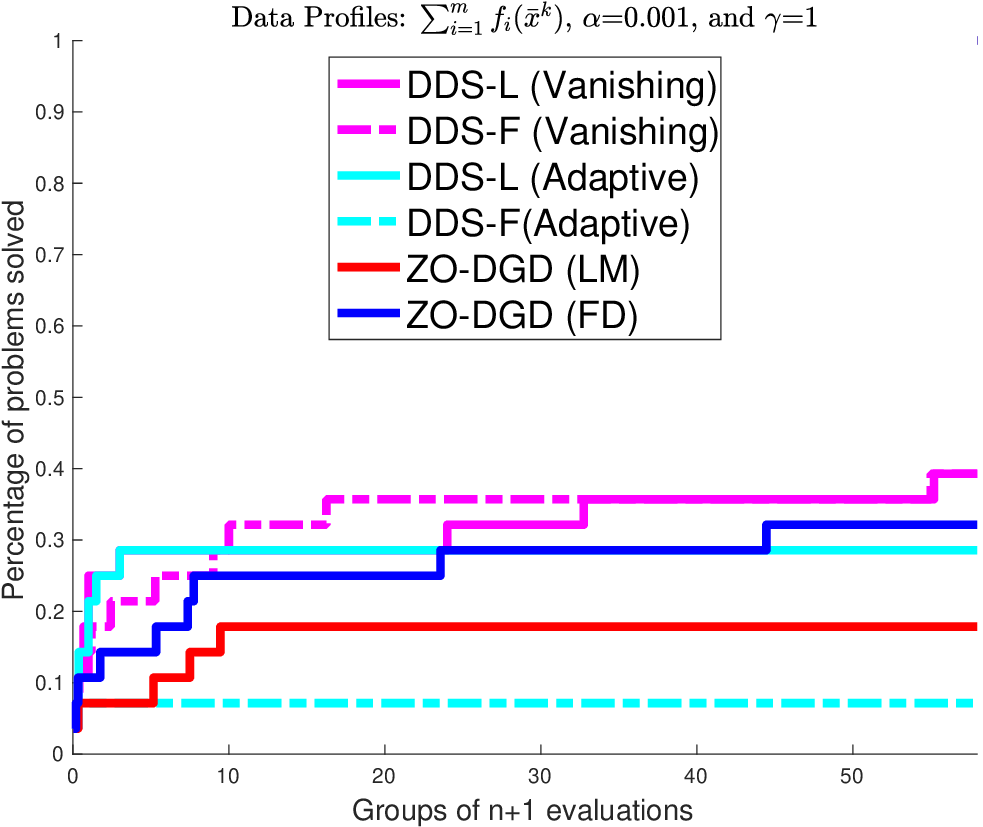}
\includegraphics[scale=0.3]{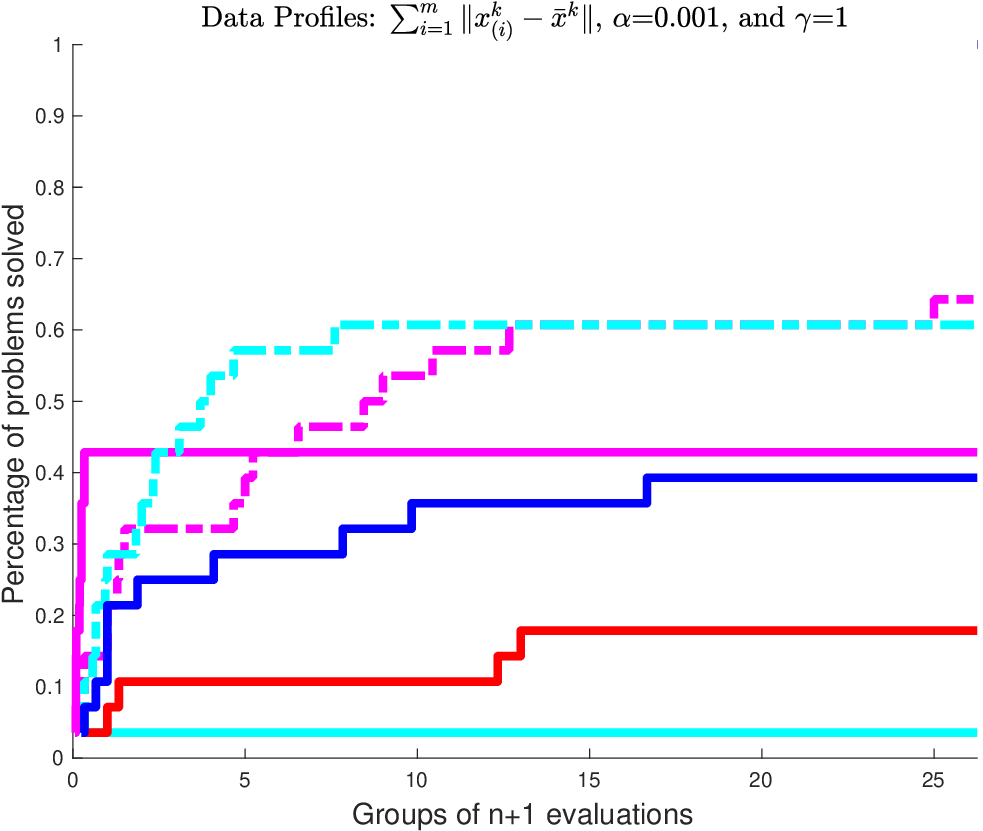}
}
\subfigure[$\alpha=10^{-6}$ \label{dp:01:1e-06}]{
\includegraphics[scale=0.3]{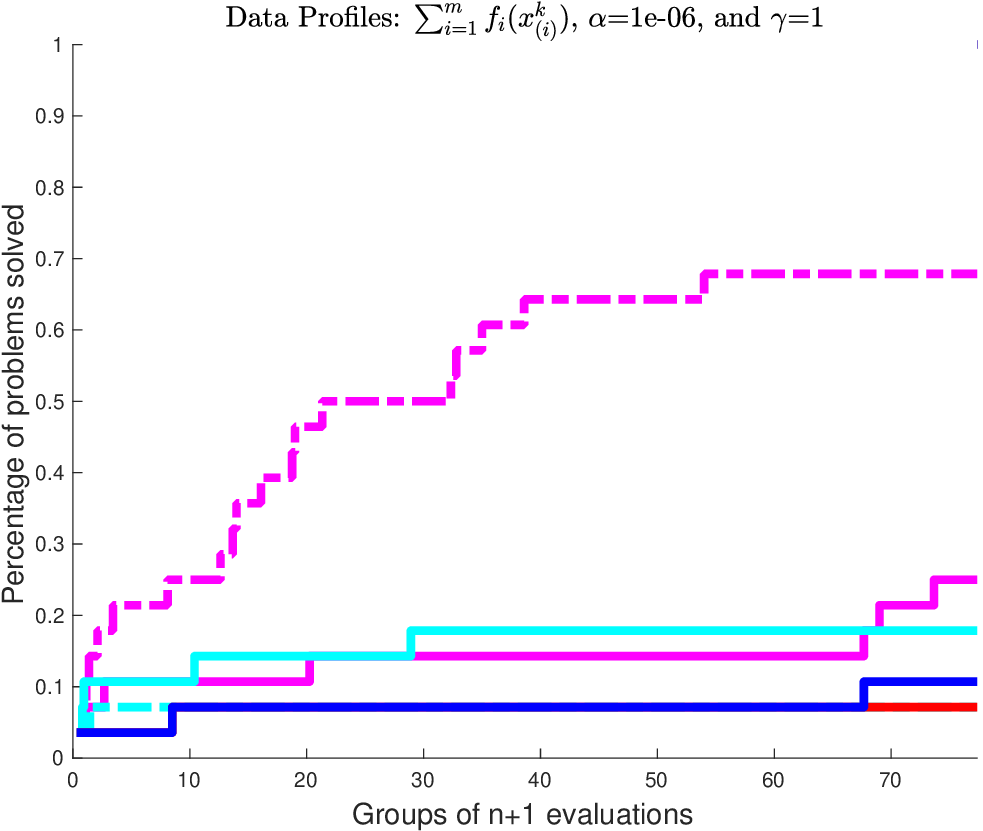}
\includegraphics[scale=0.3]{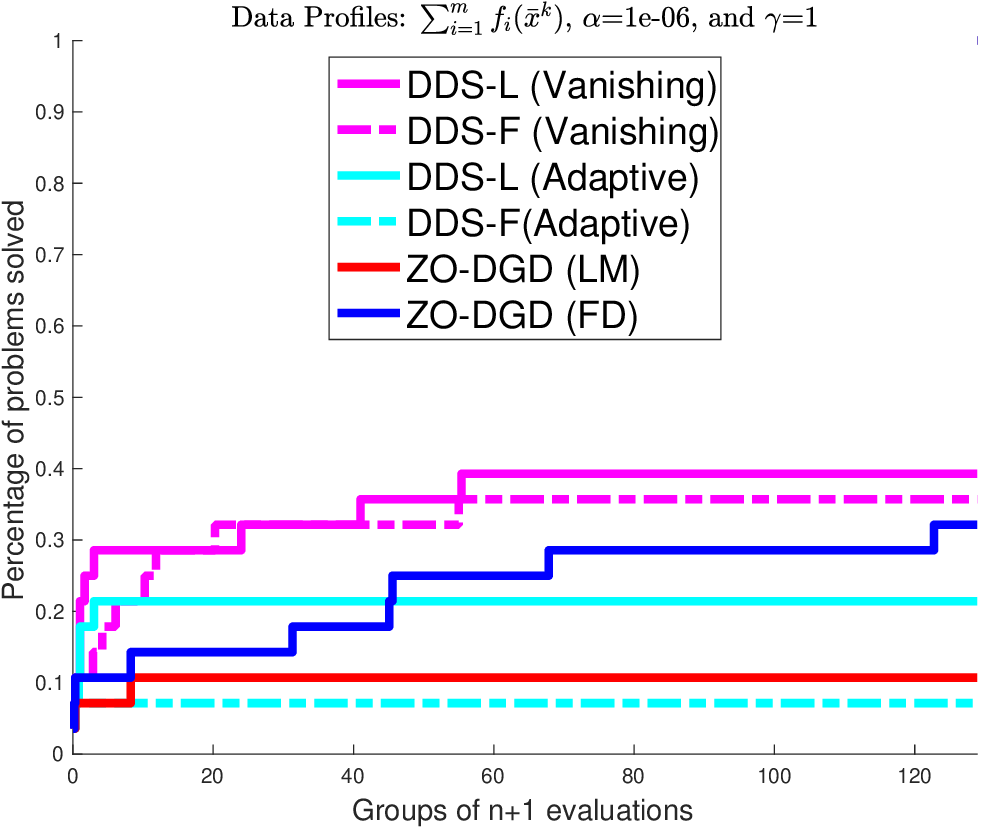}
\includegraphics[scale=0.3]{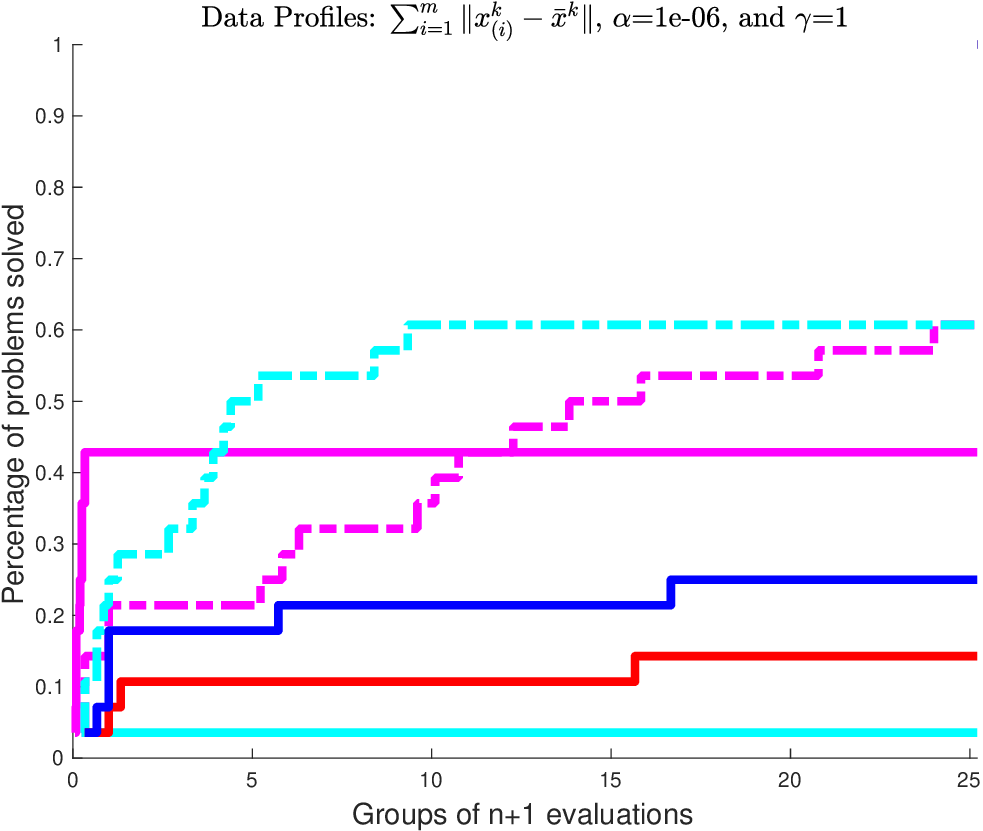}
}
\caption{Data profiles  using three different optimality metrics.}
\label{dp:fig:1e-3}
\end{figure}

Figure~\ref{dp:fig:1e-3} complements our study by presenting data 
profiles~\cite{JJMore_SMWild_2009} for our runs. Those profiles are 
consistent with the performance profiles, and further illustrate that 
\texttt{DDS-F (Vanishing)} reaches the best compromise between function 
value and consensus. Overall, these experiments support the use of 
direct-search techniques in a decentralized setting.

\section{Conclusion}
\label{sec:conc}

In this paper, we adapted direct-search techniques to operate in a decentralized 
setting. We proposed sufficient decrease conditions and stepsize updating 
techniques that borrow from the decentralized gradient descent literature as 
well as the derivative-free optimization literature. While only endowed with 
partial convergence guaranties, our algorithms can outperform zeroth-order 
decentralized gradient descent techniques in practice

Our study can be extended in several directions. First, other decentralized schemes, such as gradient tracking algorithms, could be combined with 
direct-search techniques. Extending our framework to account for
nonsmoothness or stochasticity in the objective values is also an interesting
avenue for future research. 

\bibliographystyle{plain}
\bibliography{refs}

\begin{thebibliography}{10}

\bibitem{AuHa2017}
C.~Audet and W.~Hare.
\newblock {\em {Derivative-Free and Blackbox Optimization}}.
\newblock Springer Series in Operations Research and Financial Engineering.
  Springer, Cham, Switzerland, 2017.

\bibitem{conn2009introduction}
A.R. Conn, K.~Scheinberg, and L.N. Vicente.
\newblock {\em {Introduction to Derivative-Free Optimization}}.
\newblock MOS-SIAM Series on Optimization. SIAM, Philadelphia, 2009.

\bibitem{dang2024adaptive}
Q.~Dang, S.~Yang, Q.~Liu, and J.~Ruan.
\newblock Adaptive and communication-efficient zeroth-order optimization for
  distributed internet of things.
\newblock {\em IEEE Internet of Things Journal}, 11(22):37200--37213, 2024.

\bibitem{di2016next}
P.~Di~Lorenzo and G.~Scutari.
\newblock {NEXT}: {I}n-network nonconvex optimization.
\newblock {\em IEEE Trans. Signal Inform. Process. Netw.}, 2(2):120--136, 2016.

\bibitem{EDDolan_JJMore_2002}
E.~D. {Dolan} and J.~J. {Mor\'{e}}.
\newblock Benchmarking optimization software with performance profiles.
\newblock {\em Math. Program.}, 91:201--213, 2002.

\bibitem{KJDzahini_FRinaldi_CWRoyer_DZeffiro_2024}
K.~J. {Dzahini}, F.~Rinaldi, C.~W. {Royer}, and D.~Zeffiro.
\newblock Revisiting theoretical guarantees of direct-search methods.
\newblock arXiv:2403.05322v2, 2024.

\bibitem{gao2023decentralized}
H.~Gao, M.~T. Thai, and J.~Wu.
\newblock When decentralized optimization meets federated learning.
\newblock {\em IEEE Network}, 37(5):233--239, 2023.

\bibitem{ghadimi2013stochastic}
S.~Ghadimi and G.~Lan.
\newblock Stochastic first- and zeroth-order methods for nonconvex stochastic
  programming.
\newblock {\em SIAM J. Optim.}, 23(4):2341--2368, 2013.

\bibitem{Gratton_2015}
S.~Gratton, C.~W. Royer, L.~N. Vicente, and Z.~Zhang.
\newblock Direct search based on probabilistic descent.
\newblock {\em SIAM J. Optim.}, 25(3):1515--1541, 2015.

\bibitem{JDGriffin_TGKolda_RMLewis_2008}
J.~D. {Griffin}, T.~G. {Kolda}, and R.~M. {Lewis}.
\newblock Asynchronous parallel generating set search for linearly constrained
  optimization.
\newblock {\em SIAM J. Sci. Comput.}, 30:1892--1924, 2008.

\bibitem{hajinezhad2017zeroth}
D.~Hajinezhad, M.~Hong, and A.~Garcia.
\newblock {ZONE: Z}eroth order nonconvex multi-agent optimization over
  networks.
\newblock {\em IEEE Trans. Automat. Control}, 64:3995--4010, 2019.

\bibitem{PDHough_TGKolda_VTorczon_2001}
P.D. {Hough}, T.~G. {Kolda}, and V.~Torczon.
\newblock Asynchronous parallel pattern search for nonlinear optimization.
\newblock {\em SIAM J. Sci. Comput.}, 23:134--156, 2001.

\bibitem{TGKolda_2005}
T.~G. {Kolda}.
\newblock Revisiting asynchronous parallel pattern search for nonlinear
  optimization.
\newblock {\em SIAM J. Optim.}, 16:563--586, 2005.

\bibitem{kolda2003directsearch}
T.~G. {Kolda}, R.~M. {Lewis}, and V.~Torczon.
\newblock Optimization by direct search: {N}ew perspectives on some classical
  and modern methods.
\newblock {\em SIAM Rev.}, 45:385--482, 2003.

\bibitem{JLarson_MMenickelly_SMWild_2019}
J.~Larson, M.~Menickelly, and S.~M. {Wild}.
\newblock Derivative-free optimization methods.
\newblock {\em Acta Numer.}, 28:287--404, 2019.

\bibitem{li2021communication}
Z.~Li and L.~Chen.
\newblock Communication-efficient decentralized zeroth-order method on
  heterogeneous data.
\newblock In {\em 2021 13th International Conference on Wireless Communications
  and Signal Processing (WCSP)}, pages 1--6. IEEE, 2021.

\bibitem{liu2024decentralized}
Y.~Liu, T.~Lin, A.~Koloskova, and S.~U. {Stich}.
\newblock Decentralized gradient tracking with local steps.
\newblock {\em Optim. Methods Softw.}, pages 1--28, 2024.

\bibitem{JJMore_SMWild_2009}
J.~J. Mor{\'e} and S.~M. Wild.
\newblock Benchmarking derivative-free optimization algorithms.
\newblock {\em SIAM J. Optim.}, 20:172--191, 2009.

\bibitem{nedic2009distributed}
A.~Nedic and A.~Ozdaglar.
\newblock Distributed subgradient methods for multi-agent optimization.
\newblock {\em IEEE Trans. Automat. Control}, 54(1):48, 2009.

\bibitem{EKRyu_WYin_2022}
E.~K. {Ryu} and W.~Yin.
\newblock {\em Large-Scale Convex Optimization: Algorithms \& Analyses via
  Monotone Operators}.
\newblock Cambridge University Press, 2022.

\bibitem{sahu2018distributed}
A.~K. Sahu, D.~Jakovetic, D.~Bajovic, and S.~Kar.
\newblock Distributed zeroth order optimization over random networks: {A
  Kiefer-Wolfowitz} stochastic approximation approach.
\newblock In {\em 2018 IEEE Conference on Decision and Control (CDC)}, pages
  4951--4958. IEEE, 2018.

\bibitem{sahu2020decentralized}
A.~K. Sahu and S.~Kar.
\newblock Decentralized zeroth-order constrained stochastic optimization
  algorithms: Frank--wolfe and variants with applications to black-box
  adversarial attacks.
\newblock {\em Proceedings of the IEEE}, 108(11):1890--1905, 2020.

\bibitem{shah2024adaptive}
S.~M. {Shah}, Al.~S. {Berahas}, and R.~Bollapragada.
\newblock Adaptive consensus: {A} network pruning approach for decentralized
  optimization.
\newblock {\em SIAM J. Optim.}, 34(4):3653--3680, 2024.

\bibitem{sun2016distributed}
Y.~Sun, G.~Scutari, and D.~Palomar.
\newblock Distributed nonconvex multiagent optimization over time-varying
  networks.
\newblock In {\em 2016 50th Asilomar Conference on Signals, Systems and
  Computers}, pages 788--794. IEEE, 2016.

\bibitem{tang2020distributed}
Y.~Tang, J.~Zhang, and N.~Li.
\newblock Distributed zero-order algorithms for nonconvex multiagent
  optimization.
\newblock {\em IEEE Trans. Control Netw. Syst.}, 8(1):269--281, 2020.

\bibitem{vicente2013complexity}
L.~N. {Vicente}.
\newblock Worst case complexity of direct search.
\newblock {\em EURO J. Comput. Optim.}, 1:143--153, 2013.

\bibitem{yuan2016convergence}
K.~Yuan, Q.~Ling, and W.~Yin.
\newblock On the convergence of decentralized gradient descent.
\newblock {\em SIAM J. Optim.}, 26(3):1835--1854, 2016.

\bibitem{zeng2018nonconvex}
J.~Zeng and W.~Yin.
\newblock On nonconvex decentralized gradient descent.
\newblock {\em IEEE Trans. Signal Process.}, 66(11):2834--2848, 2018.

\end{thebibliography}

\end{document}